%% file: main.tex
\documentclass[letterpaper, 10 pt, journal, twoside]{IEEEtran}

\IEEEoverridecommandlockouts                            

\usepackage[pdftex]{graphicx}
\graphicspath{{figure/}}
\usepackage{amsmath}
\usepackage{amssymb}
\usepackage{url}

\usepackage{amsthm} 
\usepackage{cite}
\usepackage{color}
\usepackage{comment}
\usepackage{fancyhdr} 
\usepackage{multirow}

\usepackage{algorithm}
\usepackage{algorithmic}
\usepackage{bm}
\usepackage{bbm}
\usepackage{enumitem}
\usepackage{mathtools}

\newtheorem{theorem}{\textbf{Theorem}}

\newtheorem{remark}{\textbf{Remark}}

\newcommand{\bx}{\bm{x}}

\newcommand{\rG}{\mathrm{G}}
\newcommand{\rL}{\mathrm{L}}
\newcommand{\rl}{\mathrm{loss}}
\newcommand{\rN}{\mathrm{N}}
\newcommand{\rB}{\mathrm{B}}
\newcommand{\rC}{\mathrm{C}}

\newcommand{\by}{\bm{y}}
\newcommand{\bz}{\bm{z}}

\newcommand{\bw}{\bm{w}}

\newcommand{\R}{\mathbb{R}}

\newcommand{\sN}{\mathcal{N}}
\newcommand{\sX}{\mathcal{X}}
\newcommand{\sE}{\mathcal{E}}

\newcommand{\sG}{\mathcal{G}}

\newcommand{\sL}{\mathcal{L}}

\newcommand{\sJ}{\mathcal{J}}

\newcommand{\sT}{\mathcal{T}}

\newcommand{\coo }{\ensuremath{\mathrm{CO_2}}}

\allowdisplaybreaks 

\title{Carbon-Aware Optimal Power Flow}

\author{Xin Chen, Andy Sun, Wenbo Shi, Na Li
\thanks{X. Chen is with the Department of Electrical and Computer Engineering, Texas A\&M University, USA; correspondence email: {xin\_chen@tamu.edu}. A. Sun is with the Sloan School of Management and the MIT Energy Initiative, Massachusetts Institute of Technology, USA. W. Shi is with Singularity Energy Inc., USA. N. Li is  with the School of Engineering and Applied Sciences, Harvard University, USA.}
\thanks{ 
 This work was supported by the NSF ASCENT Award No. 2328241, the Texas A\&M Blockchain \& Energy Research Consortium, the Texas A\&M Smart Grid Center, and the MIT Future Energy Systems Center.} 
}

\begin{document}

\maketitle


	
\input{Abstract}

\input{Nomenclature}

\input{Introduction}

\input{Preliminary}

\input{Problem}

\input{ES}

\input{Simulation}

\input{Conclusion}

\input{Appendix}

\input{Reference}

\input{BioNoPhoto}

\end{document}

%% file: Abstract.tex
\begin{abstract}
To facilitate effective decarbonization of the electric energy sector, this paper introduces a generic Carbon-aware Optimal Power Flow (C-OPF) methodology for power system decision-making that considers the active management of the grid's carbon footprints. Built upon conventional Optimal Power Flow (OPF) models, the proposed C-OPF model further integrates carbon emission flow equations and constraints, as well as carbon-related objectives, to co-optimize electric power flow and carbon emission flow across the power grid. Essentially, the proposed C-OPF can be viewed as a carbon-aware generalization of OPF. Moreover, this paper rigorously establishes the conditions that guarantee the feasibility and solution uniqueness of the carbon emission flow equations, and it proposes a reformulation technique to address the critical issue of undetermined power flow directions in the C-OPF model. Furthermore, two novel carbon footprint models for energy storage systems are developed and incorporated into the C-OPF method. Numerical simulations demonstrate the characteristics and effectiveness of the C-OPF method, in comparison with conventional OPF solutions.

\end{abstract}


%% file: Nomenclature.tex
\section*{Nomenclature}

\addcontentsline{toc}{section}{Nomenclature}

\subsection{Sets and Parameters}
\begin{IEEEdescription}[\IEEEusemathlabelsep\IEEEsetlabelwidth{$V_1,V_2,V_3$}]
\item [$\sN_i$] Set of neighbor nodes of node $i$.
\item [$\sN_i^+ (\sN_i^-)$] Set of neighbor nodes that send power to (receive power from) node $i$.
\item [$\sG_i$] Set of generators at node $i$.
\item [$\sL_i$] Set of loads at node $i$.
\item [$\sT$] Set of time steps with the time interval $\delta_t$.
\item [$w_{i,g}^\mathrm{G}$] Generation carbon 
emission factor of generator $g$ at node $i$. 
\end{IEEEdescription}	

\subsection{Variables}
\begin{IEEEdescription}[\IEEEusemathlabelsep\IEEEsetlabelwidth{$V_1,V_2,V_3$}]
\item [$R_{ij}^i (R_{ij}^j)$] Carbon flow rate of branch $ij$ from node $i$ to node $j$ measured at node $i$ (node $j$). 
\item [$R_{ij}^{\mathrm{loss}}$] Carbon flow rate associated with the power loss of branch $ij$. 
\item [$R_{i,g}^\mathrm{G}$] Carbon flow rate from generator $g$ at node $i$. 
\item [$R_{i,l}^\mathrm{L}$] Carbon flow rate injected to load $l$ at node $i$. 
\item [$P_{ij}^i (P_{ij}^j)$] Active power flow of branch $ij$ from node $i$ to node $j$ measured at node $i$ (node $j$). 
\item [$P_{ij}^{\mathrm{loss}}$] Power loss of branch $ij$. 
\item [$P_{i,g}^\mathrm{G} (Q_{i,g}^\mathrm{G})$] Active (reactive) power output of generator $g$ at node $i$.
\item [$P_{i,l}^\mathrm{L} (Q_{i,l}^\mathrm{L})$] Active (reactive) power of load $l$ at node $i$.
\item [$P_{i,t}^{\mathrm{ch}} (P_{i,t}^{\mathrm{dc}})$] Charging (discharging) power of the ES system at node $i$ at time $t$.
\item [$e_{i,t}^{\mathrm{es}}$] Energy of  the ES system at node $i$ at time $t$.
\item [$E_{i,t}^{\mathrm{es}}$] Virtually stored carbon emissions of  the ES system at node $i$ at time $t$.
\item [$w_{i,t}^{\mathrm{es}}$] Internal carbon emission intensity of  the ES system at node $i$ at time $t$.
\item [$w_i$] Nodal carbon intensity of node $i$.
\end{IEEEdescription}
\noindent
\textbf{Notes}: 1) Notations with an additional subscript $t$ denote the values at time $t$. For example, $w_{i,t}$ denotes the nodal carbon intensity of node $i$ at time $t$. 2) For a matrix $\bm{A}$, $\bm{A}[i,j]$ denotes the element in $i$-th row and $j$-th column.

%% file: Introduction.tex
\section{Introduction} \label{sec:introduction}

\IEEEPARstart{D}{eep} and rapid decarbonization of electric power systems has emerged as an urgent priority \cite{portner2022climate} to combat climate change. In 2022, the U.S. electric power sector emitted 1,539 million tons of carbon dioxide (\coo), which accounts for over 30\% of the total U.S. energy-related carbon emissions \cite{eia2023}. To enable transparent and effective grid decarbonization, precisely  
measuring and quantifying the amount of carbon emissions (i.e., carbon footprints) associated with electricity production and consumption, known as \emph{carbon accounting} \cite{world2014scope2}, is crucial. It lays the quantitative foundation necessary for informing decarbonization decisions, carbon-electricity markets,  regulation and policy development. Although almost all carbon emissions in power systems physically originate from 
electric generators due to the combustion of fossil fuels, 
it is the electricity consumption that creates the need for power generation and results in emissions. Hence, in addition to accurately measuring the emissions produced by electric generators,  it is essential to determine the carbon footprints of end-users by 
appropriately {attribute} the generation-side emissions to end-users based on their electricity use, which is referred to as \emph{demand-side carbon accounting} \cite{ourvisionpaper}. 
Accordingly, the Greenhouse Gas (GHG) Protocol\footnote{The GHG Protocol \cite{world2004greenhouse,world2014scope2} developed by the World Resources Institute (WRI)  provides internationally recognized GHG accounting and reporting standards and guidelines, which are widely used in the industry.}  
\cite{world2004greenhouse,world2014scope2} establishes two categories of emissions, i.e., Scope 1 and Scope 2, to distinguish the direct generation-side emissions and the indirect (or attributed)  emissions associated with electricity consumption and power network losses.

As introduced in \cite{ourvisionpaper}, carbon accounting frameworks can be categorized into two primary types: \emph{attributional} and \emph{consequential}. \emph{Attributional} carbon accounting aims to allocate direct generation-side emissions to end-users for assigning emissions responsibility, namely the Scope 2 emissions described above. In contrast, \emph{consequential} carbon accounting seeks to 
evaluate the change or impact on grid emissions resulting from specific decisions or projects, compared to a counterfactual baseline emission
scenario, which employs methods such as \emph{marginal emissions} \cite{10049684} and avoided emissions. These two carbon accounting frameworks 
are distinct and serve different purposes \cite{ourvisionpaper}. In terms of attributional carbon accounting,  
there are two main currently-used methods \cite{world2014scope2}: 1) the location-based method, which adopts grid \emph{average} emission factors (AEFs) across long time horizons and large areas \cite{totalCA2021,world2014scope2} to account for electricity users' carbon footprints, and 2) the market-based method, which derives carbon emission factors purely based on clean power market instruments \cite{brander2018creative}, such as  Renewable Energy Certificates (RECs) \cite{lau2008bottom} and Power Purchase Agreement (PPA) \cite{kansal2018introduction}. 
See \cite{ourvisionpaper} for a
comprehensive overview of carbon accounting methods for power grids.
 
 This paper focuses on the location-based method for attributional demand-side carbon accounting, while the existing AEF-based schemes suffer from two critical limitations: 1)  {lack of temporal and spatial granularity}, and 2) {disregard of actual electricity delivery through physical power networks} \cite{ourvisionpaper}. Since the generation fuel mix in power systems constantly changes over time, the grid carbon intensities of electricity are dynamic and time-varying with significant daily and seasonal patterns \cite{miller2022hourly}. In addition, the grid carbon intensities vary geographically, as generators of different fuel types are distributed across various locations. Moreover, power grids feature specific network topologies and circuits through which power flows, physically connecting end-users with generators and impacting carbon emissions. Therefore, carbon accounting schemes require sufficient granularity and alignment with physical power grids to 1) reflect the temporal variations and spatial diversity in grid emissions, 2) provide precise and accurate carbon accounting results for end-users, and 3) effectively inform and incentivize grid decarbonization decisions. See \cite{ourvisionpaper,totalCA2021,247epri2022} for more discussions and justifications.

To tackle these issues, the concept of \emph{carbon (emission) flow} is introduced in \cite{kang2012carbon}, where carbon emissions are treated as \emph{virtual} network flows embodied in energy flows, transmitted from producers to consumers. References \cite{kang2015carbon,li2013carbon} establish the mathematical models of carbon flow in electric power networks. The carbon flow method defines nodal and branch carbon intensities for the grid, 
providing a temporally and spatially granular depiction of the grid's carbon footprints. In this way, 
the carbon flow method represents a prominent tool for demand-side carbon accounting, which aligns carbon footprint calculations with the physical grids and power flows. See   \cite{kang2015carbon} and Section \ref{sec:preliminary} for more details.  Recent work \cite{chen2024contributions} 
presents a tree search algorithm to trace the contribution of each generator to individual lines and loads, thereby estimating nodal carbon
emissions. This approach builds on the concept of electricity flow tracing \cite{bialek1996tracing}  and yields nodal carbon emission estimates comparable to those produced by the carbon flow method, with the potential distinction being in how power losses are handled.

Unlike the studies above 
 that focus on accounting for and estimating the grid carbon footprints, this paper aims to advance foundational methodologies for grid decarbonization decision-making and enable optimal carbon-electricity joint management in power grids. In this regard, the \emph{Optimal Power Flow (OPF)} method \cite{frank2012optimal,cain2012history} stands as a foundational mathematical tool for optimizing power system decisions. It has been studied extensively in the literature and widely applied in grid planning, operations, control, and electricity markets \cite{chen2019aggregate,chen2020distributed,chen2017robust}. Existing OPF schemes typically seek optimal power grid decisions to minimize a specific economic cost objective, while satisfying the power flow equations,  network constraints, and operational limits of devices. However, the essential goals of grid decarbonization and carbon footprint management are \emph{inadequately} considered. Since power system decisions, e.g., the siting and sizing of renewable generators or power dispatch of various generation sources, 
directly impact the grid's carbon footprints, 
it becomes necessary to explicitly integrate carbon footprint management into grid decision-making for achieving desired decarbonization performances and outcomes.

\subsection{Related Work and Key Issues}

There have been
a number of recent studies \cite{sun2017analysis,shen2020low,cheng2018bi,wu2022carbon,wei2023wasserstein,gu2022carbon,wang2021optimal,sang2023encoding,wang2022robust,lu2022peer} that consider carbon emission flow in 
power system planning and operation.
 Reference \cite{sun2017analysis} proposes a transmission expansion planning method that defines an index to quantify the equity performance of carbon emission allocation based on the carbon flow model. 
In \cite{shen2020low}, a multi-objective power network transition model is built to plan the 
retirement of aging coal-fired power plants, while one of the objectives is to minimize user-side carbon footprints.  Reference \cite{8804208}  studies the low-carbon operation of multiple energy systems and derives the locational energy-carbon integrated price based on the nodal carbon intensities calculated using the carbon flow method. 
In \cite{cheng2018bi,wu2022carbon,wei2023wasserstein,gu2022carbon}, carbon-aware expansion planning models are established for multi-energy systems under carbon emission constraints on electric devices and energy hubs. Additionally, the carbon flow model and constraints are taken into account in power scheduling
\cite{wang2021optimal}, energy management \cite{sang2023encoding,wang2022robust}, and peer-to-peer carbon-electricity trading \cite{lu2022peer}. 
Existing studies outlined above focus on specific power system applications. However, it remains unaddressed in establishing a fundamental decision-making methodology necessary for guiding various grid decarbonization decisions and supporting theoretical studies and performance analysis.

Moreover, there are two critical issues in the integration of  carbon emission flow into the grid decision models that necessitate to be addressed:
\begin{itemize}
    \item [1)] (\emph{Power Flow Directions}): The carbon flow model \cite{kang2015carbon} needs to pre-determine the power flow directions for all branches to identify the power inflows for each node. However,  the directions of branch power flows are typically unknown prior to solving optimal decision models.  Most existing works that employ the carbon flow model in grid decision-making either overlook the issue of unknown power flow directions or 
    assume they are predetermined by alternating between grid optimization and carbon flow calculation through iterations \cite{lu2022peer}. 
    Reference \cite{gu2022carbon} introduces binary indicator variables to handle the unknown power flow directions in the carbon flow model. It results in a mixed-integer nonconvex quadratically constrained optimization problem, which cannot be readily solved using off-the-shelf optimizers due to involving both integer variables and nonconvex constraints. Thus, a tailored heuristic penalty-based iterative algorithm is designed to solve the optimization problem.

\item[2)] (\emph{Carbon Footprint Models for Energy Storage}): By switching between charging and discharging, energy storage (ES) systems can shift load demand and transfer renewable energy across time, offering substantial potential to reduce power system emissions. 
Hence, developing carbon footprint models for ES systems is crucial, since it lays the quantitative foundation for making optimal ES planning and operation decisions, such as determining installation location and capacity, and managing charging and discharging. In addition, granular and accurate carbon accounting for ES systems enhances information transparency,  supporting the development of low-carbon regulatory policies and market mechanisms for ES. It also incentivizes ES stakeholders to adopt low-carbon practices to reduce their own carbon footprint and overall grid emissions.
Carbon accounting for ES systems has attracted considerable recent attention from the industry \cite{esaccount2021,esacc2022}. References \cite{wang2021optimal,gu2022carbon,gu2023carbon}  propose different carbon footprint models for ES systems, 
while they neglect the carbon emission leakage associated with ES energy loss (see Remark \ref{remark:ES}), and the carbon accounting mechanisms for ES system owners remain unclear. 
\end{itemize}

\subsection{Our Contributions}

In this paper, we introduce a generic \emph{Carbon-aware Optimal Power Flow (C-OPF)} methodology as the fundamental theory for carbon-aware power system decision-making that incorporates the active management of the
grid’s carbon footprints. Built upon conventional OPF models, 
the C-OPF model (see model \eqref{eq:copf}) further integrates the carbon flow equations and constraints, as well
as carbon-related objectives, to
co-optimize power flow
and carbon flow in the power grid. Essentially, C-OPF is a carbon-aware generalization of OPF, and it produces optimal decisions that satisfy carbon emission constraints and balance the power-related and carbon-related costs. 
The main contributions of this paper are threefold:  
\begin{itemize}
    \item [1)] To our knowledge, this is the first work that introduces the generic C-OPF methodology and builds its mathematical model for lossy power networks with formulation examples. In particular, this paper rigorously establishes the conditions that guarantee the feasibility and solution uniqueness of the carbon flow equations (see Theorem \ref{thm:inverse}), and presents the key properties of the C-OPF model. 
 
    \item [2)] We propose a reformulation approach to address the issue of unknown power flow directions in the C-OPF model by introducing dual power flow variables
    and complementarity constraints. This reformulation is exact and eliminates the necessity of knowing power flow directions in advance.

    \item [3)] We develop two novel carbon footprint models for ES systems: one precisely models the dynamics of carbon emissions virtually stored in an ES unit and the carbon leakage associated with ES energy loss, while the other treats ES as a load during charging and a carbon-free generator during discharging. The corresponding carbon accounting schemes for ES owners are also introduced.

\end{itemize}

Furthermore, we develop a carbon-aware economic dispatch model as an example based on the proposed C-OPF method, and 
demonstrate the effectiveness of C-OPF through numerical experiments in comparison with OPF-based solutions.

\begin{remark} (Key Merits of C-OPF). \normalfont
 In contrast to conventional OPF-based schemes that merely incorporate carbon emission costs into the objective,  our proposed C-OPF method possesses three distinct merits: 1) C-OPF explicitly models carbon flow alongside power flow,  enabling a \emph{granular} representation and management of the grid's carbon footprints rather than focusing solely on the system-wide total emissions. 2) C-OPF is flexible to model various global and local decarbonization targets or regulatory requirements for different entities and stakeholders. It can ensure that power system decisions comply with these requirements by imposing corresponding carbon flow constraints.  
 3) C-OPF  integrates demand-side carbon accounting mechanisms that attribute emissions from power generation to consumption, which allows the optimization of carbon footprints for end-users, e.g., via carbon-aware power dispatch and demand response. It establishes the theoretical foundation to engage numerous end-users with substantial power flexibility and resources in grid decarbonization decision-making.
  
\end{remark}

The remainder of this paper is organized as follows: 
Section \ref{sec:preliminary} introduces the concept and model of carbon flow as well as the use for carbon accounting. Section \ref{sec:problem} presents the C-OPF method with the reformulation approach. Section \ref{sec:ES} introduces the carbon footprint models for energy storage systems. 
  Numerical experiments are conducted in Section \ref{sec:simulation}, and conclusions are drawn in Section \ref{sec:conclusion}.

%% file: Preliminary.tex
\section{Carbon Emission Flow and Carbon Accounting}\label{sec:preliminary}

In this section, we first introduce the concept and model of carbon flow for lossy power networks, and then establish the conditions that ensure the feasibility and solution uniqueness of the carbon flow equations. Next, we present
the application of the carbon flow method for demand-side carbon accounting.

\subsection{Concept of Carbon Emission Flow}   \label{sec:cfconcept}

\begin{figure}
      \centering
      \includegraphics[scale=0.44]{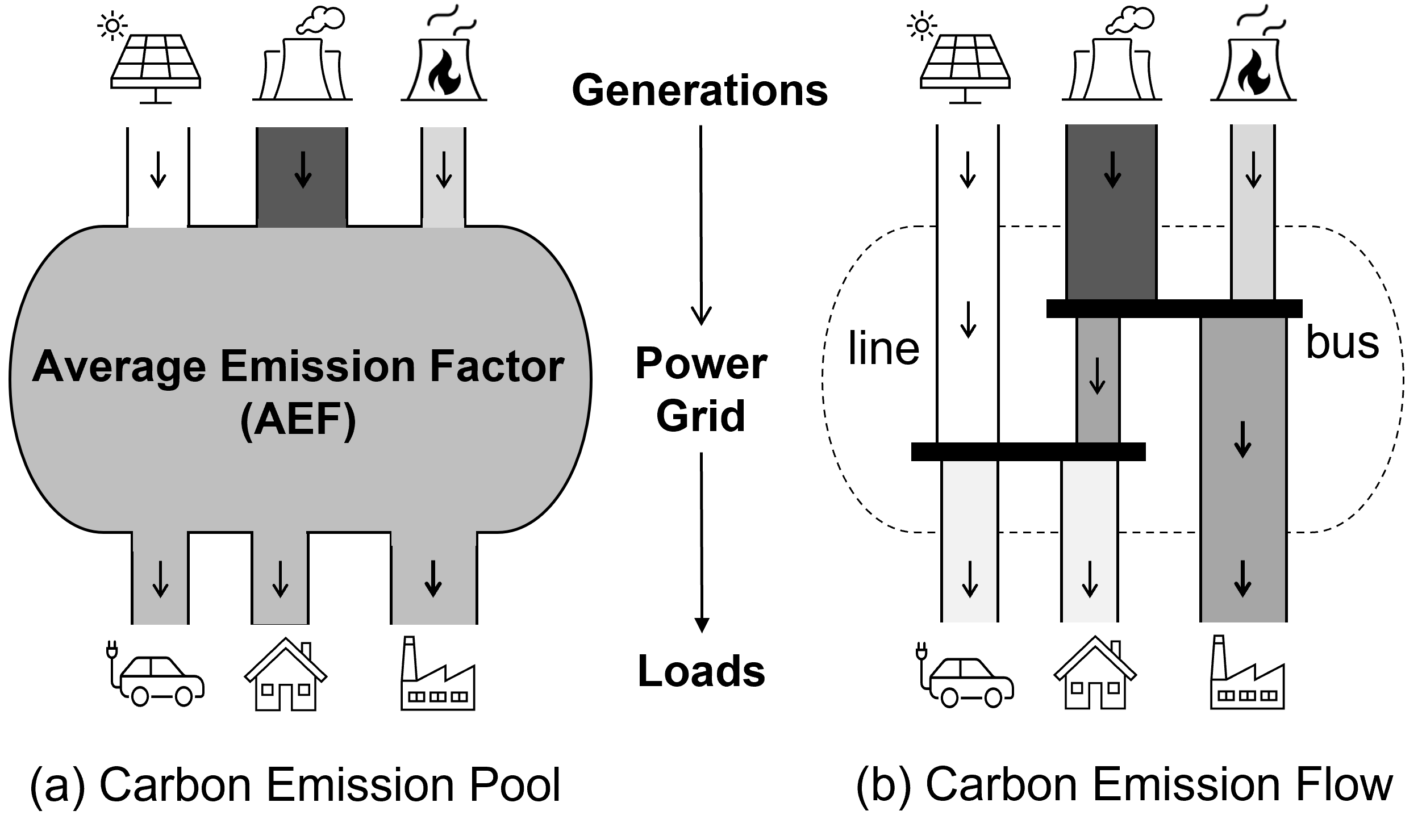}
      \caption{Comparison between carbon emission pool and carbon emission flow. (In sub-figure (a), all end-users in a large area adopt the same grid average emission factor (AEF) to calculate their attributed carbon footprints. 
      In sub-figure (b), each pipeline represents a power line, with the width indicating the magnitude of power flow.
      Darker colors indicate higher carbon emission intensities. Power in-flows with different carbon emission intensities are mixed at each bus and distributed downstream.)}
      \label{fig:cflow}
\end{figure}

As mentioned above, the location-based method calculates the average emission factor (i.e., the ratio of total carbon emissions to the total generation energy) of an area grid over a defined period for demand-side carbon accounting. As illustrated in Figure \ref{fig:cflow}, this AEF method treats the power system as a large ``pool", assuming that all generations and end-users are connected to one homogeneous and shared infrastructure, without consideration of the power network and power flow. In contrast, the carbon flow method views the carbon emissions from generators as \emph{virtual} attachments to the power flow. These emissions are considered to be transmitted from generators through power networks and accumulate on the user side, forming carbon flows.
 The concept of carbon flow is analogous to a ``\emph{water supply}'' system, where virtual carbon emissions accompanying power flows can be viewed as invisible particles contained in water flows that are delivered from water sources to users. The carbon emission intensity of electricity is analogous to the particle concentration of water, and a renewable generator is analogous to a pure water source. 
 In this way, the carbon flow method aligns with the physical power grids and underlying power flows, enabling a temporally and spatially granular depiction of the grid's carbon footprints. 

In practice, the grid operators and utilities are expected to implement the carbon flow scheme, as these entities oversee power grid operations and possess detailed power flow profiles. Moreover, the carbon flow information (e.g., real-time nodal carbon intensities) represents useful grid emission signals, which can be
sent to end-users (e.g., by displaying on smart meters) to inform user-side decarbonization decisions \cite{chen2024enhance}.

\subsection{Carbon Emission Flow Model for Lossy Power Networks} \label{sec:cfmodel}

\begin{figure}
    \centering
    \includegraphics[scale=0.33]{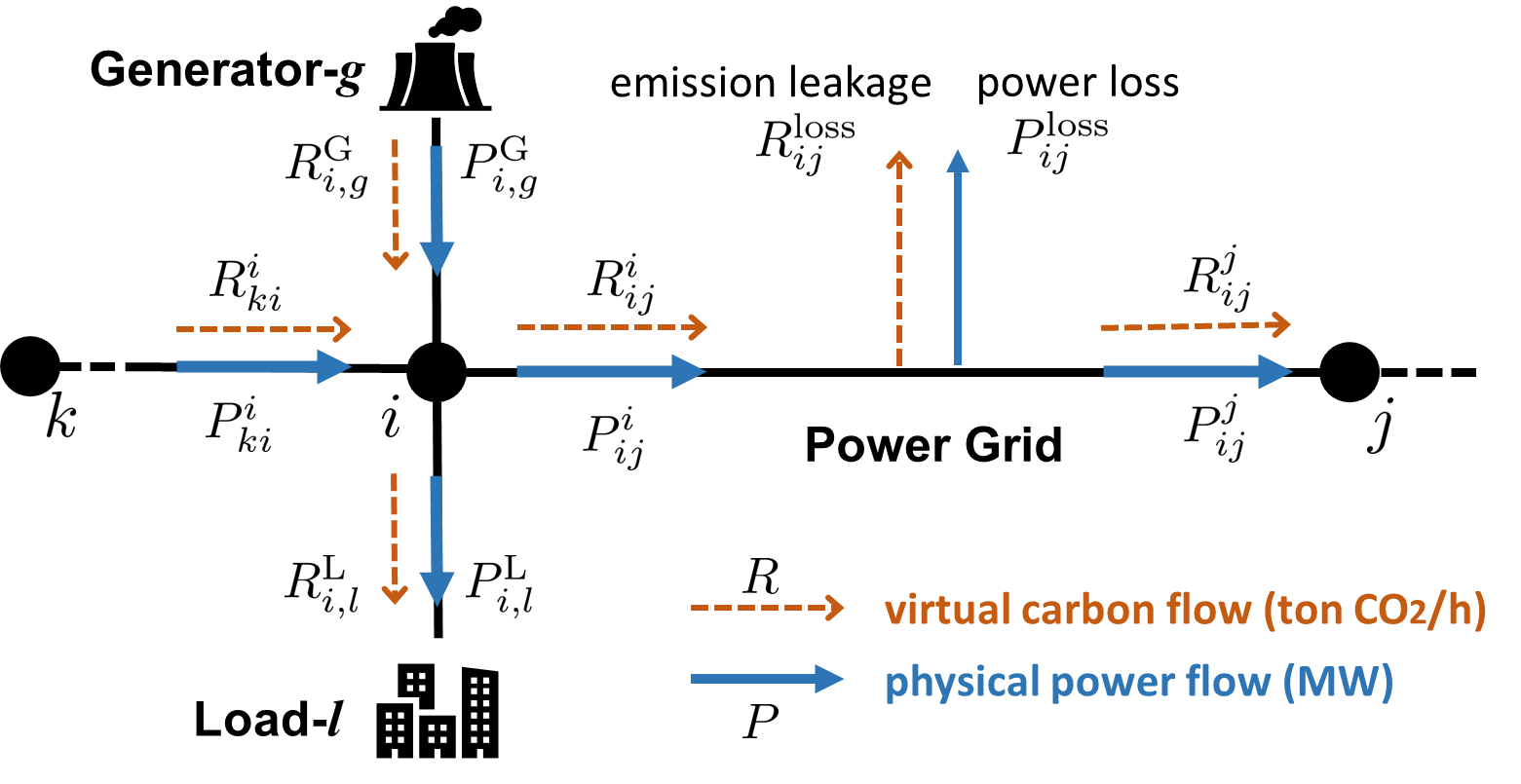}
    \caption{Illustration of virtual carbon flow and physical power flow.}
    \label{fig:cfvalue}
    \vspace{-5pt}
\end{figure}

Consider a power network described by a connected graph $G(\mathcal{N},\sE)$, where $\sN\!:=\!\{1,\cdots, N\}$ denotes the set of nodes and $\sE\subset \sN\times\sN$ denotes the set of branches. As illustrated in Figure \ref{fig:cfvalue}, the notations with $P$ denote the active power flow values (in the unit of MW), while the notations with $R$ represent the associated carbon flow rates (in the unit of ton\coo/h). 
The \emph{carbon (emission) intensity} is defined as the ratio $w \!=\! \frac{R}{P}$ (in the unit of ton\coo/MWh) that describes the amount of carbon flow associated with one unit of power flow. 

 The active power flows adhere to the power balance equation \eqref{eq:powerbalance} at each node $i$:
\begin{align} \label{eq:powerbalance}
      \sum_{k\in \mathcal{N}_i^+}\!\! P_{ki}^i + \sum_{g\in \mathcal{G}_i}\! P_{i,g}^\mathrm{G} =     \sum_{j\in \mathcal{N}_i^-} \!\! P_{ij}^i + \sum_{l\in \mathcal{L}_i}\! P_{i,l}^{\rL},&&  \forall i\in\sN, 
\end{align}
and $P_{ij}^i =  P_{ij}^{\mathrm{loss}}+P_{ij}^j $ for each branch $ij\in \sE$. Then, the carbon flow model is built upon the active power flows and the following two fundamental principles \cite{kang2015carbon}:
\begin{itemize}
    \item [1)] \emph{Conservation Law of Nodal Carbon Mass}: Similar to the power flow balance \eqref{eq:powerbalance},  the total carbon inflows equal the total carbon  outflows at each node $i\in \sN$, i.e., 
    \begin{align}\label{eq:conservation}
      \sum_{k\in \mathcal{N}_i^+}\!\! R_{ki}^i + \sum_{g\in \mathcal{G}_i}\! R_{i,g}^\mathrm{G} =     \sum_{j\in \mathcal{N}_i^-} \!\! R_{ij}^i + \sum_{l\in \mathcal{L}_i}\! R_{i,l}^{\rL}.
    \end{align}
    and $R_{ij}^i = R_{ij}^{\rl} + R_{ij}^j$ for each branch $ij\in \sE$. In \eqref{eq:conservation}, $R_{i,g}^{\rG}\!=\! w_{i,g}^{\rG} P_{i,g}^{\rG}$ is the  generation carbon emission rate of generator-$g$ at node $i$.  
    \item [2)] \emph{Proportional Sharing Principle}: 
    At each node $i\in\sN$,
     the allocation of total carbon inflows among all outflows is proportional to their active power flow values, i.e., 
 \begin{align}
      R_{ij}^i = \frac{R_i^{\mathrm{in}}}{P_i^{\mathrm{in}}}\cdot P_{ij}^i, \ \  R_{i,l}^{\rL}   =  \frac{R_i^{\mathrm{in}}}{P_i^{\mathrm{in}}}\cdot P_{i,l}^{\rL},
 \end{align}
where $P_i^{\mathrm{in}}\!:= \!  \sum_{k\in \mathcal{N}_i^+}\! P_{ki}^i +\sum_{g\in \mathcal{G}_i}\! P_{i,g}^\mathrm{G}  $ 
and  $R_i^{\mathrm{in}}\!:=\!  \sum_{k\in \mathcal{N}_i^+}\! R_{ki}^i +  \sum_{g\in \mathcal{G}_i}\! R_{i,g}^\mathrm{G} $
denote the total power inflow and total carbon inflow at node $i\in\sN$.
\end{itemize}
 The \emph{nodal carbon intensity} of each node $i$ is calculated as \eqref{eq:mci}:
\begin{align}\label{eq:mci}
  w_i  \!= \! \frac{R_i^{\mathrm{in}}}{P_i^{\mathrm{in}}}\!=\!  \frac{ \sum_{g\in\sG_i} w_{i,g}^{\rG} P_{i,g}^{\rG} + \sum_{k\in \sN_i^+ } w_k P_{ki}^i}{\sum_{g\in\sG_i}  P_{i,g}^{\rG} + \sum_{k\in \sN_i^+ } P_{ki}^i },\ \forall i\in\sN,
\end{align}
which is the ratio of total carbon inflow to the total power inflow.
Equation \eqref{eq:mci} indicates that the power flow and power loss of each branch $ij\in\sE$ share the same carbon emission intensity that equals to the nodal carbon intensity of the sending node, i.e.,  $ \frac{R_{ij}^i}{P_{ij}^i} =  \frac{R_{ij}^j}{P_{ij}^j} =  \frac{R_{ij}^{\rl}}{P_{ij}^{\rl}} =w_i$.

Equation \eqref{eq:mci} is referred to as the \emph{carbon flow model} or \emph{carbon flow equation}, and it can be equivalently reformulated in the matrix form \eqref{eq:mci:ma}. 
 The detailed derivation of \eqref{eq:mci:ma1} is provided in 
Appendix \ref{app:cfder}. 
\begin{subequations}\label{eq:mci:ma}
    \begin{align}
 & (\bm{P}_{\rN} - \bm{P}_{\rB}  ) \bm{w}_{\rN} =   \bm{r}_{\rG} \label{eq:mci:ma1}\\
&\qquad  \Longrightarrow\   \bm{w}_{\rN} = (\bm{P}_{\rN} - \bm{P}_{\rB}  )^{-1}  \bm{r}_{\rG}, \label{eq:mci:ma2}
\end{align}
\end{subequations}
Here, $\bm{w}_{\rN}\!:=\! (w_i)_{i\in\sN}\!\in\!\R^{N}$ and $\bm{r}_{\rG} \!:=\!(\sum_{g\in \mathcal{G}_i}\! \! w_{i,g}^{\rG} P_{i,g}^{\rG})_{i\in\sN}\!\in\!\R^{N} $ 
denote the column vectors that collect the nodal carbon intensities and nodal generation carbon flows, respectively.
$\bm{P}_{\rN}\!:=\!\text{diag}(P_i^{\mathrm{in}})\!\in\! \R^{N\times N} $ is the diagonal matrix whose $i$-th diagonal entry is the nodal active power inflow $P_i^{\mathrm{in}}$ to node $i$. $\bm{P}_{\rB} \!\in\! \R^{N\times N}$ is the branch power inflow matrix that is built by letting $\bm{P}_{\rB}[i,k] \!=\! P_{ki}^i$ and $\bm{P}_{\rB}[k,i] \!= \!0$ if node $k$ sends power flow $P_{ki}^i$ to node $i$. 
Note that $P_{ki}^i$ and $P_{ki}^k$ are different due to the line power loss. 
Given power flow profiles,
the linear equations \eqref{eq:mci:ma1} can be solved to obtain $\bm{w}_{\rN}$, e.g., through matrix factorization and backward and forward substitution, rather than directly computing the matrix inverse in \eqref{eq:mci:ma2}.

\subsection{Feasibility and Solution Uniqueness of Carbon Flow Model}

A critical fundamental question is whether the carbon flow equations \eqref{eq:mci} or \eqref{eq:mci:ma} for a power network are feasible and admit a unique solution. Based on the matrix form \eqref{eq:mci:ma}, it translates to the invertibility of the carbon flow matrix $\bm{P}_{\rC}:=\bm{P}_{\rN} - \bm{P}_{\rB}$, which implies the feasibility and solution uniqueness. Hence, this subsection studies the properties of the matrix $\bm{P}_{\rC}$.

By definition, the matrix $\bm{P}_{\rC}$ is diagonally dominant \cite{golub2013matrix}, and we define the set $\mathcal{J}$:
\begin{align} \label{eq:diagdom}
    \mathcal{J} := \Big\{i\in \sN\, \Big\rvert\, |\bm{P}_{\rC}[i,i]| > \!\!\sum_{j=1;j\neq i}^N\! |\bm{P}_{\rC}[i,j]| \Big\}.
\end{align}
In practice, $\mathcal{J}\neq \emptyset$, i.e., there exist some rows  $i\in\sJ$ of $\bm{P}_{\rC}$ that are strictly diagonally dominant; such rows correspond to the nodes with generation power injections $P_{i,g}^{\rG}$. 
Then, we establish the \emph{invertibility} property of $\bm{P}_{\rC}$ in Theorem \ref{thm:inverse}.

\begin{theorem}\label{thm:inverse}
Suppose that for each $i\notin \sJ$ of the matrix $\bm{P}_{\rC}$, there is a sequence of nonzero elements of $\bm{P}_{\rC}$ of the form $\bm{P}_{\rC}[i,i_1], \bm{P}_{\rC}[i_1,i_2],\cdots, \bm{P}_{\rC}[i_r, j]$ with $j\in \sJ$. Then, 
 $\bm{P}_{\rC}$ is invertible, and $\bm{P}_{\rC}$ is an {M-matrix}.
\end{theorem}
\begin{proof}
According to \cite[Theorem]{shivakumar1974sufficient}, the matrix $\bm{P}_{\rC}$ is nonsingular and thus invertible. Since $\bm{P}_{\rC}$ is diagonally dominant and all its off-diagonal elements $\{\bm{P}_{\rC}[i,j]\}_{i\neq j}$ are non-positive,  $\bm{P}_{\rC}$ is an \emph{M-matrix} \cite{ando1980inequalities} according to \cite[Corollary 4]{shivakumar1974sufficient}. 
\end{proof}

Theorem \ref{thm:inverse} indicates that under the condition in Theorem \ref{thm:inverse}, the matrix $\bm{P}_{\rC}$ is invertible, and thus
the carbon flow equations are feasible and have a unique solution regarding the nodal carbon intensities and other carbon flow values. 

\begin{figure}
    \centering
    \includegraphics[scale=0.42]{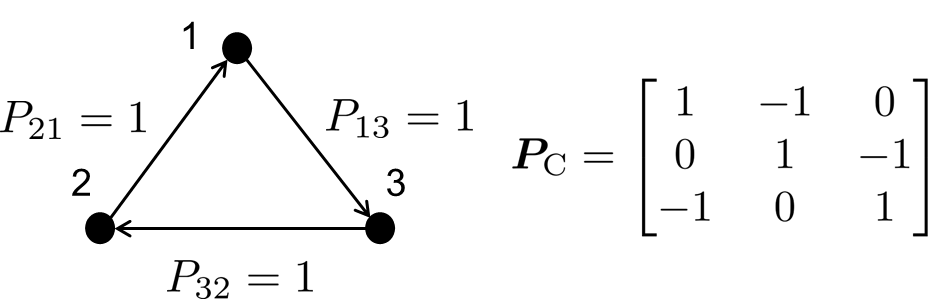}
    \caption{A simple 3-node network case example.}
    \label{fig:3bus}
\end{figure}

\begin{remark}
\normalfont
    The condition in Theorem \ref{thm:inverse} can be interpreted as that for any node $i$ with no generation power injection, one can find a power flow path $i \leftarrow i_1 \leftarrow i_2 \leftarrow \cdots \leftarrow i_r \leftarrow j$ in the power network that traces upstream to a node $j$ that has generation power injection. This condition
{generally} holds for practical connected power networks. This condition also implies that every node has non-zero power flux, i.e., $\bm{P}_{\rC}[i,i] >0$ for all $i\in[N]$. However, non-zero nodal power flux can \emph{not} guarantee that $\bm{P}_{\rC}$ is invertible. 
For example, consider a simple 3-node circular network case shown in Figure \ref{fig:3bus}. In this case, every node has 1 unit of power flowing through it, but the matrix $\bm{P}_{\rC}$ is singular. And this case violates the condition in Theorem \ref{thm:inverse}. In addition, $\bm{P}_{\rC}$ possesses all the properties of being an \emph{M-matrix}. For example, all the elements of $\bm{P}_{\rC}^{-1}$ are non-negative and every eigenvalue of $\bm{P}_{\rC}$ has a positive real part \cite[Theorem 1.1]{ando1980inequalities}. 
\end{remark}

 \subsection{Use of Carbon Emission Flow for Carbon Accounting}\label{sec:accounting}

Given power flow results, one can
solve the set of carbon flow equations \eqref{eq:mci} or \eqref{eq:mci:ma} to obtain the nodal carbon intensities $\bm{w}_{\rN}$, and then calculate all carbon flow values. According to the GHG Protocol \cite{world2004greenhouse,world2014scope2}, the carbon accounting rules are: 
{\setlist{leftmargin=4mm}
\begin{itemize}
    \item   Generator-$g$ at node $i$ shall account for the (Scope 1) direct carbon emission rate $R_{i,g}^{\rG} = w_{i,g}^{\rG} P_{i,g}^{\rG}$;
    \item   Load-$l$ at node $i$ shall account for the (Scope 2) attributed carbon emission rate  $R_{i,l}^{\rL} = w_i P_{i,l}^{\rL}$; 
    \item   The power network owner shall account for the (Scope 2) attributed carbon emission rate associated with the network power loss, i.e., $\sum_{ij\in\sE}\! R_{ij}^{\rl} = \sum_{ij\in\sE}\! w_i P_{ij}^{\rl}$.
\end{itemize}}

As the generation fuel mix and power flows change over time, the carbon flow profiles are time-varying, and 
the carbon footprints over a time period $\sT$ (e.g., one day or one year) are the time accumulation of $R$. For example, the carbon footprint of load-$l$ at node $i$ is calculated as  $\hat{E}_{i,l,\sT}:=\sum_{t\in \sT} \delta_t w_{i,t}P_{i,l,t}^{\rL}$, where $\delta_t$ is the time interval (e.g., 1 hour or 15 minutes).

%% file: Problem.tex
\section{Carbon-Aware Optimal Power Flow Method} \label{sec:problem}

In this section, we introduce the generic C-OPF model
with examples, and present the reformulation technique to tackle the power flow direction issue of C-OPF.

\subsection{Generic  C-OPF Model}

To enable the optimal management of carbon footprints in power system decision-making, we propose the generic C-OPF model \eqref{eq:copf} to co-optimize power flow and carbon flow.
\begin{subequations} \label{eq:copf}
\begin{align}
\text{Obj.}  \ \ &  \min_{\bm{x} \in \mathcal{X}} \  f_{\mathrm{power}}(\bm{x},\bm{y}) + f_{\mathrm{carbon}} (\bx,\by,\bm{z}) \label{eq:copf:obj}\\
 \text{s.t.} \ \   & \text{Power Flow Equations }(\bm{x}, \bm{y}) = \bm{0},\label{eq:copf:pf}\\
    & \text{Power Flow Constraints } (\bm{y})\leq \bm{0}, \label{eq:copf:pfc}\\
        & \text{Carbon Flow Equations } (\bm{x}, \bm{y}, \bm{z}) = \bm{0},\label{eq:copf:cf} \\
    & \text{Carbon Flow Constraints } (\bx, \by, \bm{z})\leq \bm{0}. \label{eq:copf:cfc}
\end{align}
\end{subequations}
Here, $\bx$ denotes the decision variables subject to the feasible set $\mathcal{X}$. 
The specification of $\bx$ and $\sX$ depends on the practical applications. For instance, 
$\bx$ denotes the generation decisions of various generators and $\sX$ represents the generation capacity limits and ramping constraints in economic dispatch \cite{476062}.
$\bx$ can also be the load adjustment decisions in demand response, or the site and size decisions of new renewable generators in grid planning. 
$\by$ denotes the power flow-related variables, such as network voltage profiles and branch power flows. $\bz$ represents the carbon flow-related variables, such as nodal carbon intensities $\bm{w}_{\rN}$ and carbon flow values.

 \subsubsection{Objective Function} 
 The objective \eqref{eq:copf:obj} aims to minimize the overall cost that consists of two components: the power-related cost denoted as $f_{\mathrm{power}}$ and the carbon emission-related cost $f_{\mathrm{carbon}}$. 
  Depending on the specific applications,  $f_{\mathrm{power}}$ can be the generation cost, network loss, grid expansion investment cost, etc. $f_{\mathrm{carbon}}$ is defined to capture the externality of carbon emissions and regulatory penalty on generation-side or demand-side emissions (or the bonus on emission reduction).
 An example of these cost functions is given by \eqref{eq:obj}:
 \begin{subequations} \label{eq:obj}
       \begin{align} 
      f_{\mathrm{power}} & := \sum_{i\in\sN} \sum_{g\in \sG_i} \Big(c_{i,g}^2  (P_{i,g}^{\rG})^2 + c_{i,g}^1  {P_{i,g}^{\rG}} + c_{i,g}^0\Big), \label{eq:cpower}\\
      f_{\mathrm{carbon}} &:= c^{\mathrm{emi}}\cdot \sum_{i\in\sN} \sum_{g\in \sG_i} {w_{i,g}^{\rG}} {P_{i,g}^{\rG}},\label{eq:ccarbon}
  \end{align}
 \end{subequations}
  where \eqref{eq:cpower} denotes the total  generation cost in a quadratic form with the parameters $c_{i,g}^2, c_{i,g}^1, c_{i,g}^0$,
  and \eqref{eq:ccarbon} is the penalty on generation-side emissions with the cost coefficient $c^{\mathrm{emi}}$.

 \subsubsection{Power Flow Equations and Constraints}  
The power flow equations
 \eqref{eq:copf:pf} and power flow constraints \eqref{eq:copf:pfc} remain the same as them in   classic OPF models \cite{giannakis2013monitoring,cain2012history,frank2012optimal}. The full AC power flow equations
are formulated as \eqref{eq:powerflow:full}:
 \begin{subequations} \label{eq:powerflow:full}
     \begin{align}
       &  P_{ij}^i     =  (V_i^2 -V_iV_j\cos(\theta_i-\theta_j))g_{ij} \nonumber\\
         &\qquad\qquad\qquad -V_iV_j\sin(\theta_i-\theta_j)b_{ij},&& \forall ij\in \sE,   \label{eq:pf:P}
\\
 & Q_{ij}^i    =   (V_iV_j\cos(\theta_i-\theta_j)-V_i^2)b_{ij}\nonumber \\
         &\qquad\qquad\qquad  -V_iV_j \sin(\theta_i-\theta_j)g_{ij},&&  \forall ij\in \sE, \\
 & \sum_{j\in \sN_i} P_{ij}^i\, = \sum_{g\in \sG_i} P_{i,g}^{\rG} - \sum_{l\in \sL_i} P_{i,l}^{\rL},&&  \forall i\in\sN,  \label{eq:full:p} \\
    &   \sum_{j\in \sN_i} Q_{ij}^i = \sum_{g\in \sG_i} Q_{i,g}^{\rG} - \sum_{l\in \sL_i} Q_{i,l}^{\rL},&&   \forall i\in\sN,
     \end{align}
 \end{subequations}
 where $V_i$ and $\theta_i$ are the voltage magnitude and phase angle at node $i$. $Q_{ij}^i$ denotes the reactive power flow of branch $ij$ from node $i$ to node $j$ measured at node $i$. $g_{ij}$ and $b_{ij}$ are the  conductance and susceptance of branch $ij$.

 The power flow constraints \eqref{eq:copf:pfc} generally involve 
the line thermal capacity constraints \eqref{eq:thermal} and voltage limits \eqref{eq:voltage}.
\begin{subequations} \label{eq:pf:con}
    \begin{align} 
    \qquad    & (P_{ij}^i)^2 + (Q_{ij}^i)^2\leq \Bar{S}_{ij}^2,&& \forall ij\in \sE, \label{eq:thermal} \\
       \qquad           &\quad \underline{V}_i \leq \, V_i \leq \Bar{V}_i,&& \forall i\in\sN,  \label{eq:voltage}
    \end{align}
\end{subequations}
where $\Bar{S}_{ij}$ is the apparent power capacity of line $ij$. $\Bar{V}_i, \underline{V}_i$ are the upper and lower limits of voltage magnitude at node $i$.

 Since the power flow equations and constraints remain unchanged, the linearization and convexification methods developed for them are still applicable. For instance, the classic DC power flow model \eqref{eq:dc} \cite{giannakis2013monitoring}, which neglects power loss, can be used in the C-OPF model to replace the complete power flow equations \eqref{eq:powerflow:full} for simplicity. In \eqref{eq:dc:pf2}, the superscript ``$i$" is omitted since
 $P_{ij}^i = P_{ij}^j$ in the DC power flow model due to neglecting power loss.
\begin{subequations} \label{eq:dc}
\begin{align}
&P_{ij}   = -b_{ij}\cdot(\theta_i - \theta_j), && \forall ij\in\sE, \label{eq:dc:pf2} \\
 & \text{Equation }\eqref{eq:full:p}.
\end{align}
\end{subequations}

\subsubsection{Carbon Flow Equations and Constraints} 
 The carbon flow equations \eqref{eq:copf:cf} are given by \eqref{eq:mcieq}:
\begin{align}\label{eq:mcieq}
     w_i &\big(\sum_{g\in\sG_i}\!  P_{i,g}^{\rG} +  \sum_{j\in \sN_i^+}\!\!  P_{ji}^i  \big) \nonumber \\
     &\qquad=\! \sum_{g\in\sG_i}\! w_{i,g}^{\rG} P_{i,g}^{\rG} +\!\! \sum_{j\in \sN_i^+}\!\!\! w_j P_{ji}^i,\qquad \forall i \in\sN,
\end{align}
which is simply an equivalent reformulation of \eqref{eq:mci}.

The carbon flow constraints \eqref{eq:copf:cfc} can take various forms depending on the practical settings for carbon footprint management. For example, an upper limit $\bar{\bw}_{\rN}:=(\bar{w}_i)_{i\in\sN}$ can be imposed on nodal carbon intensities, i.e. \eqref{eq:capnci}, to ensure that users at these nodes are supplied with low-carbon electricity.
By adjusting this upper limit $\bar{\bw}_{\rN}$, one can control the level of ``cleanness" of the supplied electricity at a certain location. 
In particular, letting $\bar{w}_i = 0$ enforces that the electricity supply at node $i$ is completely carbon-free. The definition of ``node" is flexible in terms of geographical scales, which can represent a district grid, a distribution feeder, or a balancing area.
\begin{align} 
 \qquad\qquad    w_i & \leq \bar{w}_i, &&\forall i \in\sN. \label{eq:capnci}
\end{align} 
An alternative carbon flow constraint \eqref{eq:copf:cfc} can impose a cap $\bar{E}_{i,l}^{\rL}$ on the total individual user-side emissions as \eqref{eq:capne}:
\begin{align}
   \qquad \delta_t\sum_{t\in\sT} ( w_{i,t}\cdot P_{i,l,t}^{\rL}) &\leq \bar{E}_{i,l}^{\rL}, &&\forall l\in\sL_i, i\in\sN.\label{eq:capne}
\end{align}
Besides, instead of an emission cap on individual users, a cap $\bar{E}_{i}^{\rL}$ on the total nodal level emissions can be imposed as \eqref{eq:cap:node}: 
\begin{align} \label{eq:cap:node}
 \qquad  \delta_t \sum_{t\in\sT}\big(w_{i,t} \cdot\!\sum_{l\in\sL_i} P_{i,l,t}^{\rL}\big) &\leq \bar{E}_{i}^{\rL}, &&\forall i\in\sN.
\end{align} 
The determination of the cap parameters depends on practical requirements. 
Other carbon flow constraints, 
including requirement for emission allocation fairness and equity \cite{sun2017analysis}, can also be employed.

 Due to the proportional sharing principle used in the carbon flow model, a natural limit on nodal carbon intensities is $w_i\in[0,w^{\rG}_{\max}]$ for all $i\in\sN$, where  $w_{\max}^{\rG}:= \max_{i,g}\{ w_{i,g}^{\rG} \}$ is the largest generation carbon emission factor. Hence, if    $\bar{w}_i$ in  \eqref{eq:capnci} is set to be larger than $w_{\max}^{\rG}$  for some nodes $i$,   \eqref{eq:capnci} imposes no actual constraint on these nodes.
 Another critical issue is that inappropriately designed carbon flow constraints may render the  C-OPF model \eqref{eq:copf} infeasible, e.g., when the caps in \eqref{eq:capnci}-\eqref{eq:cap:node} are too small to be achievable. To address this issue, the hard constraints \eqref{eq:capnci}-\eqref{eq:cap:node} can be converted to soft constraints by adding slack variables with corresponding penalties in the objective. For example, one can replace constraint \eqref{eq:cap:node} with  \eqref{eq:cap:relax} and add a penalty term $\sum_{i\in\sN}(c_i^{\mathrm{E}}\cdot\alpha_i) $ to the carbon-related cost function $f_{\mathrm{carbon}}$ in objective \eqref{eq:copf:obj}:
 \begin{align} \label{eq:cap:relax}
   \delta_t \sum_{t\in\sT}\big(w_{i,t}  \!\sum_{l\in\sL_i} P_{i,l,t}^{\rL}\big) &\leq \bar{E}_{i}^{\rL} + \alpha_i,   \ \alpha_i\geq 0, &&\forall i\in\sN,
 \end{align}
 where $\alpha_i$ is the slack variable and $c_i^{\mathrm{E}}$ denotes the penalty cost coefficient for excessive demand-side emissions.  
 
\begin{remark} 
\normalfont 
We note that the carbon flow constraints are not inherent physical limitations but rather regulatory requirements imposed on the power grid and end-users to manage their carbon footprints. The objective is to incentivize low-carbon electricity supply and consumption behaviors and ensure compliance with grid decarbonization regulations and targets. For instance, to meet the emission cap constraints \eqref{eq:capne}, \eqref{eq:cap:node},  
end-users can optimally schedule their load trajectories $(P_{i,l,t}^{\rL})_{t\in\sT}$ to increase (decrease) electricity consumption when the grid exhibits low (high) nodal carbon intensity. Alternatively, additional renewable generation units can be deployed in proximity to node $i$, to effectively reduce the nodal carbon intensity $w_{i,t}$. These lead to the development of optimal carbon-aware demand response and expansion planning schemes based on the C-OPF method, and it can also be used to support many other carbon-aware decision applications in power systems.
 \qed
\end{remark}

 Essentially,
the C-OPF model  \eqref{eq:copf} is a carbon-aware generalization of the OPF model,  and it reduces to an OPF model if the carbon-related objective function $f_{\mathrm{carbon}}$, carbon flow equations \eqref{eq:copf:cf} and constraints \eqref{eq:copf:cfc} are removed or inactive. It also implies that existing OPF techniques, such as linearization, convexification, decomposition, stochastic modeling, etc., can still be applied to the power flow components in C-OPF \eqref{eq:copf}. 
Moreover, the  C-OPF model \eqref{eq:copf} can be directly extended to the multi-period settings and involve time-coupled constraints, such as generator ramping limits and the state-of-charge constraints of ES systems.
As an illustrative example, 
a multi-period C-OPF-based economic dispatch model \eqref{eq:ed} is established in Section \ref{sec:ed}, which includes these time-coupled constraints and employs a full AC power flow model.

In the C-OPF model \eqref{eq:copf}, the carbon flow method is used for demand-side carbon accounting to  {align the grid's carbon footprint quantification with physical power system operation and actual power flows}. 
Nevertheless, the C-OPF method is flexible and can adapt to other carbon accounting approaches by replacing the carbon flow equation \eqref{eq:copf:cf} with other valid carbon accounting mathematical models.

\subsection{Reformulation  for Power Flow Directions}
\label{sec:keyissue}

As mentioned in the introduction section and indicated by the set $\sN_i^+$, the carbon flow equation $\eqref{eq:mcieq}$ requires the pre-determination of branch power flow directions to identify the power inflows for each node. However, the directions of branch power flows are generally unknown prior to solving the C-OPF problem. To address this issue, we introduce two non-negative power flow variables $\hat{P}_{ji}^i\!\geq\! 0$ and $\hat{P}_{ij}^i\!\geq\! 0$ for each branch $ij\in \sE$ with $P_{ij}^i = \hat{P}_{ij}^i - \hat{P}_{ji}^i$.
Specifically, $\hat{P}_{ji}^i$ and $\hat{P}_{ij}^i$ denote the power flow components from node $j$ to node $i$ and from node $i$ to node $j$, respectively, both of which are measured on the side of node $i$. 
Then, we can equivalently reformulate the carbon flow equation \eqref{eq:mcieq} as \eqref{eq:n1}, and also need to replace $P_{ij}^i$ with  $\hat{P}_{ij}^i - \hat{P}_{ji}^i$ 
in the power flow equations \eqref{eq:powerflow:full} (or the DC power flow model \eqref{eq:dc}) and constraints \eqref{eq:thermal}.
\begin{subequations}\label{eq:n1}
\begin{align}
           &w_i (\sum_{g\in\sG_i}\!   P_{i,g}^{\rG} + \! \sum_{j\in \sN_i}\!  \hat{P}_{ji}^i  ) \nonumber \\
          & \qquad= \sum_{g\in\sG_i}\!\! w_{i,g}^{\rG} P_{i,g}^{\rG} +\!\! \sum_{j\in \sN_i}\!\!\! w_j \hat{P}_{ji}^i,&& \forall i\in\sN,
  \label{eq:n1:nci}\\
   &\hat{P}_{ji}^i\geq 0, \    \hat{P}_{ij}^i\geq 0, && \forall ij \in \sE, \label{eq:n1:posi}
   \\
    &\hat{P}_{ji}^i\cdot \hat{P}_{ij}^i    =0,   && \forall ij \in \sE.\label{eq:n1:com} 
\end{align}
\end{subequations}
In \eqref{eq:n1:nci}, we replace  $\sN_i^+$  (the set of neighbor nodes that send power to node $i$)   by $\sN_i$ (the set of all neighbor nodes of node $i$).  
In addition, the complementarity constraint \eqref{eq:n1:com} is added to ensure that either $\hat{P}_{ji}^i$ or $\hat{P}_{ij}^i$ must be zero for each branch $ij$. Here, a useful trick to facilitate the nonlinear optimization is to replace the complementarity constraint \eqref{eq:n1:com} by the relaxed constraint \eqref{eq:relax} \cite{fletcher2004solving}, and this relaxation is exact due to \eqref{eq:n1:posi}.
\begin{align}\label{eq:relax}
  \qquad  \hat{P}_{ji}^i\cdot \hat{P}_{ij}^i    \leq 0,\qquad\forall ij \in \sE.
\end{align} 

Alternatively, we can introduce a binary variable $\gamma_{ij}$ for $ij\in\sE$ and linearize the complementarity constraint \eqref{eq:n1:com} with
\begin{align} \label{eq:n2}
      \gamma_{ij}\in\{0,1\},\  \hat{P}_{ji}^i \leq (1-\gamma_{ij}) \bar{P}_{ij},\ \hat{P}_{ij}^i\leq \gamma_{ij}\bar{P}_{ij}.
\end{align}
Note that both reformulations via
\eqref{eq:n1} or \eqref{eq:n2} are equivalent to the original carbon flow equation \eqref{eq:mcieq}, 
but 
 the branch flow directions do not need to be known in advance.

  In \cite{gu2022carbon}, binary indicator variables are introduced to handle the unknown power flow directions in the carbon flow model. It results in a mixed-integer nonconvex quadratically constrained optimization problem, and a tailored penalty-based iterative algorithm is designed in \cite{gu2022carbon} to solve the optimization problem through a number of iterations. In contrast, our proposed dual power flow variables reformulation method with the complementarity constraints \eqref{eq:relax}
renders the C-OPF model a standard nonconvex optimization problem without integer variables. Therefore, the  C-OPF model can be directly solved using off-the-shelf nonlinear optimizers such as IPOPT \cite{biegler2009large}, without the need for designing ad hoc solution algorithms.

%% file: ES.tex
\section{Carbon Footprint Model for Energy Storage} 
\label{sec:ES}

 Energy storage (ES) systems play a critical role in decarbonizing power grids, as their operations can be optimized to curtail overall system emissions, e.g., charging when
the grid is clean and discharging when the grid is under high emissions.
Consider an ES system connected to node $i\!\in\!\sN$. 
 For time $t\in \sT:=\{1,2,\cdots,T\}$ with the time interval $\delta_t$, the dynamical ES power model \cite{chen2021leveraging} is formulated as \eqref{eq:ESpower}:
\begin{subequations} \label{eq:ESpower}
\begin{align}
& 0\leq  P_{i,t}^{\mathrm{ch}}  \leq   \bar{P}_{i}^{\mathrm{ch}},\ \ 0\leq  P_{i,t}^{\mathrm{dc}}  \leq  \bar{P}_{i}^{\mathrm{dc}},  \\
& P_{i,t}^{\mathrm{ch}}\cdot P_{i,t}^{\mathrm{dc}} = 0,\label{eq:ES:comple}\\
  &  e_{i,t+1}^{\mathrm{es}} = \kappa_i e_{i,t}^{\mathrm{es}} + \delta_t  \big(\eta_i^{\mathrm{ch}}  P_{i,t}^{\mathrm{ch}}   - \frac{1}{\eta_i^{\mathrm{dc}}} P_{i,t}^{\mathrm{dc}} \big ),\label{eq:ES:E} \\
  &  \underline{e}_i \leq e_{i,t}^{\mathrm{es}}  \leq \Bar{e}_i,\quad   e_{i,T+1}^{\mathrm{es}} = e_{i,1}^{\mathrm{es}},\quad  \forall i\in \sN, t\in\sT.
  \label{eq:ES:uplow}
\end{align}
\end{subequations}
Here, $\bar{P}_{i}^{\mathrm{ch}}$ and $\bar{P}_{i}^{\mathrm{dc}}$ are the   charging and discharging power capacities. 
$\eta_i^{\mathrm{ch}}\!\in\!(0,1]$ and $\eta_i^{\mathrm{dc}}\!\in\!(0,1]$ denote the charging and discharging efficiency coefficients, respectively. $\kappa_i\!\in\!(0,1]$ denotes the storage efficiency factor that models the loss of stored energy over time. $\Bar{e}_i$ and $\underline{e}_i$ are the upper and lower bounds of the energy level of the ES system. The complementarity constraint \eqref{eq:ES:comple} is used to enforce that an ES unit can not charge and discharge at the same time.

 In this section, we propose two carbon footprint models for ES systems: the ``water tank" model and the ``load/carbon-free generator" model. The associated carbon accounting mechanisms for ES  owners are presented as well.

\subsection{``Water Tank" Carbon Footprint Model}\label{sec:watertank} 

The ``{water tank}" model of ES is conceptually aligned with the analogy of a ``{water supply}" system used to explain the carbon flow model in Section \ref{sec:cfconcept}.  
Analogous to a water tank that stores both water and invisible particles, an ES system is viewed to store both 
electric energy $e_{i,t}^{\mathrm{es}}$ (in the unit of MWh) and \emph{virtual} carbon emissions $E_{i,t}^{\mathrm{es}}$ (in the unit of ton\coo). 
We then define the internal ES carbon intensity $w_{i,t}^{\mathrm{es}}$  as \eqref{eq:esci}:
\begin{align} \label{eq:esci}
    w_{i,t}^{\mathrm{es}} = \frac{E_{i,t}^{\mathrm{es}}}{e_{i,t}^{\mathrm{es}}}.
\end{align} 
Note that
$e_{i,t}^{\mathrm{es}}$ should not be zero to make $w_{i,t}^{\mathrm{es}}$ well-defined. 
This can be achieved by setting the lower energy bound $\underline{e}_i$ in \eqref{eq:ES:uplow} to be a small positive value rather than zero. It can also prevent numerical issues during the optimization process.

\begin{figure}
    \centering
    \includegraphics[scale=0.34]{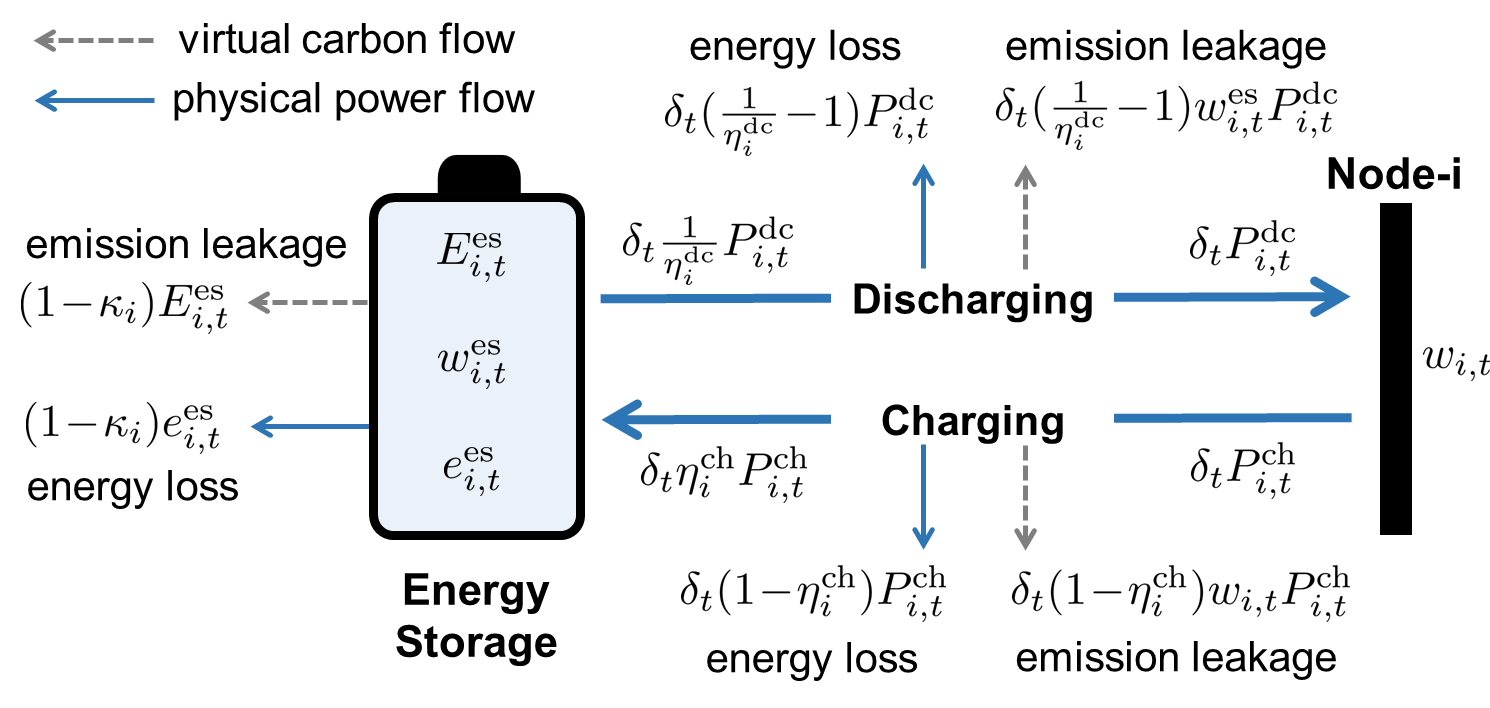}
    \caption{Power flow and carbon  flow under the ``water tank" ES model.}
    \label{fig:ES}
\end{figure}

Based on the ES power model \eqref{eq:ESpower}, we develop the  dynamical \emph{carbon footprint model of} the ES unit as \eqref{eq:EScarbon} for $t\in\sT$:
\begin{align} \label{eq:EScarbon}
    E_{i,t+1}^{\mathrm{es}}  = \kappa_i E_{i,t}^{\mathrm{es}} + \delta_t \big( w_{i,t} \eta_i^{\mathrm{ch}}   P_{i,t}^{\mathrm{ch}}  -  w_{i,t}^{\mathrm{es}}\frac{1}{\eta_i^{\mathrm{dc}}}  P_{i,t}^{\mathrm{dc}} \big).
\end{align}
Model \eqref{eq:EScarbon} implies that (virtual)  carbon emissions are injected into the ES unit when it charges with electricity in the nodal carbon intensity $w_{i,t}$; and it releases the stored emissions back to the grid when it discharges with electricity in the ES carbon intensity $w_{i,t}^{\mathrm{es}}$. In particular, \eqref{eq:EScarbon} models the \emph{carbon emission leakage} associated with the energy loss during the storage, charging, and discharging processes of the ES unit. Figure \ref{fig:ES} illustrates the 
 carbon flow and power flow of the ES unit.

Alternatively, we can plug in the ES carbon intensity \eqref{eq:esci} into \eqref{eq:EScarbon} to eliminate the variable $E_{i,t}^{\mathrm{es}}$ and reformulate the ES carbon footprint model \emph{equivalently} as \eqref{eq:ES:w}:\footnote{In \eqref{eq:ES:la}, the term $\frac{\delta_t}{    \eta_i^{\mathrm{dc}}}  P_{i,t}^{\mathrm{dc}}$ is dropped from the definition of $\lambda_{i,t}$, because 
$P_{i,t}^{\mathrm{dc}} =0$ when the ES unit charges and $\lambda_{i,t}=1$ when it discharges, which are not affected by this term. } 
\begin{subequations} \label{eq:ES:w}
      \begin{align}
  w_{i,t+1}^{\mathrm{es}}    =&   
      \lambda_{i,t} \cdot w_{i,t}^{\mathrm{es}}   +    (1- \lambda_{i,t}) \cdot w_{i,t}, \label{eq:ES:dyn2}\\
     \lambda_{i,t}  := & \frac{ \kappa_i   e_{i,t}^{\mathrm{es}} }{ \kappa_i e_{i,t}^{\mathrm{es}}  + \delta_t    \eta_i^{\mathrm{ch}} P_{i,t}^{\mathrm{ch}}  },\label{eq:ES:la}
\end{align}  
\end{subequations}
which describes the dynamics of the ES  carbon intensity $w_{i,t}^{{\mathrm{es}}}$.
Interestingly, equation \eqref{eq:ES:dyn2} indicates an intuitive feature that the ES carbon intensity $w_{i,t+1}^{\mathrm{es}}$ at the next time $t\!+\!1$ is a \emph{convex combination} of the ES carbon intensity $w_{i,t}$ at time $t$ and the nodal carbon intensity $w_{i,t}$ with the weight coefficient $\lambda_{i,t}\in [0,1]$.
When charging,  the ES unit behaves as if mixing the internally stored electricity with the newly charged electricity;  
when it discharges, $\lambda_{i,t}=1$ and  the ES carbon intensity remains unchanged, i.e., $ w_{i,t+1}^{\mathrm{es}} =  w_{i,t}^{\mathrm{es}}$.

Accordingly, the carbon flow equation \eqref{eq:n1:nci} in the C-OPF model is modified as \eqref{eq:cf:ES1} to incorporate the ES system.
\begin{align}\label{eq:cf:ES1}
     & w_{i,t} \big(  P_{i,t}^{\mathrm{dc}} +
    \sum_{g\in\sG_i}\!\!  P_{i,g,t}^{\rG} + \!\!\!\sum_{j\in \sN_i}\!\!\! \hat{P}_{ji,t}^i  \big)=  w_{i,t}^{\mathrm{es}}  P_{i,t}^{\mathrm{dc}} \nonumber\\
    &\quad  +\!
    \sum_{g\in\sG_i}\!\! w_{i,g,t}^{\rG} P_{i,g,t}^{\rG} +\!\! \sum_{j\in \sN_i}\!\!\! w_{j,t} \hat{P}_{ji,t}^i, \quad  \forall i\in\sN, t\in\sT.
\end{align} 
From \eqref{eq:cf:ES1}, it is seen that an ES unit affects the nodal carbon intensities and carbon flow only when it discharges.

 \begin{remark} \label{remark:ES}
     (Comparison with Existing ES Emission  Models). \normalfont Existing carbon emission models  
     for ES systems proposed in \cite{wang2021optimal,gu2022carbon,gu2023carbon} also formulate the virtually stored emissions and internal ES carbon intensity. However, these models 
     neglect the carbon emission leakage associated with the energy loss during the storage, charging, and discharging processes.
This issue makes these ES carbon emission models not rigorous or even problematic. For example, consider the scenario when an ES unit remains idle, i.e.,  neither charging nor discharging. 
Over time, the stored energy $e_t$ gradually depletes to zero due to energy loss, while the virtually stored carbon emissions $E_t$ remain constant as the carbon leakage associated with energy loss is not considered. As a result, the internal ES carbon intensity $E_t/e_t$  would approach an infinitely large value, which is unreasonable. In contrast, our proposed ``water tank" ES carbon footprint model \eqref{eq:EScarbon} or \eqref{eq:ES:w} avoids these issues by precisely modeling carbon leakage, ensuring that carbon footprint attribution is consistent with the actual
electric energy usage.
\qed 
 \end{remark}

\begin{remark}
    ({Carbon Accounting for ES Owners}). \normalfont Under the ``water tank" carbon footprint model, for the time horizon $\sT$, the owner of the ES system shall account for the (Scope 2) attributed carbon emission $\hat{E}_{i,\sT}^{\mathrm{es}}$  that is calculated by \eqref{eq:ES:account1}:
\begin{align}\label{eq:ES:account1}
\quad    \hat{E}_{i,\sT}^{\mathrm{es}} = \sum_{t=1}^T \delta_t(w_{i,t}P_{i,t}^{\mathrm{ch}} - w_{i,t}^{\mathrm{es}}P_{i,t}^{\mathrm{dc}}), \quad \forall i\in\sN,
\end{align} 
which is the net carbon emissions withdrawn from the grid. Intuitively, the attributed emission $ \hat{E}_{i,\sT}^{\mathrm{es}}$ can be decomposed into two parts: 1) the change of virtually stored carbon emissions, i.e., $E_{i,T+1}^{\mathrm{es}} \!-\! E_{i,1}^{\mathrm{es}}$, and 2) the carbon emission leakage associated with ES energy loss, i.e., $E_{i,1}^{\mathrm{es}} +\sum_{t=1}^T \delta_t(w_{i,t}P_{i,t}^{\mathrm{ch}} - w_{i,t}^{\mathrm{es}}P_{i,t}^{\mathrm{dc}}) - E_{i,T+1}^{\mathrm{es}} $. 
To make it more clear, we consider a lossless ES unit with $\kappa_i \!= \!\eta_i^{\mathrm{ch}}\! =\! \eta_i^{\mathrm{dc}} \!=\!1$.  
If this ES unit recovers the initially stored emission level in the final time step, namely $E_{i,T+1}^{\mathrm{es}} \!=\! E_{i,1}^{\mathrm{es}}$, we have $ \hat{E}_{i,\sT}^{\mathrm{es}} = E_{i,T+1}^{\mathrm{es}} - E_{i,1}^{\mathrm{es}} = 0$ from \eqref{eq:EScarbon}. In this case,
the ES 
owner accounts for \emph{zero} carbon emissions, regardless of the number of charging and discharging cycles. This outcome aligns with the role of an ES system, which does not directly produce or consume electricity (carbon emissions) itself but rather enables the temporal shifts.
\qed
\end{remark}

\subsection{``Load/Carbon-Free Generator" Carbon Footprint Model}

The  proposed ``water tank" model above
requires continuous monitoring of the virtually stored carbon emissions within an ES system.  To facilitate practical implementation, we propose an alternative model termed ``\emph{load/carbon-free generator}" to characterize the carbon footprints of an ES unit. Specifically, 
this model directly treats an ES unit as a load during charging and as a carbon-free clean generator during discharging. 
 Accordingly, the carbon flow equation \eqref{eq:n1:nci} in the C-OPF model is modified as \eqref{eq:cf:ES1} with $w_{i,t}^{\mathrm{es}} \equiv 0$ for all $t$, since the ES unit is regarded as a carbon-free generator during discharging.

In terms of carbon accounting, 
under the ``load/carbon-free generator" model, the ES owner shall account for the  (Scope 2) attributed carbon emissions $\hat{E}_{i,\sT}^{\mathrm{es}}$ calculated by
\eqref{eq:ES:account2}:
\begin{align}\label{eq:ES:account2}
   \hat{E}_{i,\sT}^{\mathrm{es}} = \sum_{t=1}^T \delta_tw_{i,t}P_{i,t}^{\mathrm{ch}}, \qquad \forall i\in\sN.
\end{align}
It implies that the virtual carbon emissions absorbed from the grid during the charging period
accumulate locally at the ES unit and are not released back to the grid. Thus,
the ES owner incurs higher Scope-2 carbon emissions compared with the carbon accounting scheme \eqref{eq:ES:account1} under the ``water tank" model. 
Nevertheless, under the ``load/carbon-free generator" model, ES owners can make profits by 
acting as clean energy suppliers in the carbon-electricity market, e.g., 
selling renewable energy certificates (RECs) \cite{lau2008bottom}, 
since the discharging power is considered carbon-free. As a result, ES owners are incentivized to charge their ES units when the grid is clean with low carbon intensity to reduce their own carbon footprints \eqref{eq:ES:account2}; and they are incentivized to discharge when the grid is in high emission to make more profit, as the price of clean electricity or RECs is expected to be higher at that time. In this way, the ``load/carbon-free generator" model is easy to use and naturally incentivizes the carbon-aware operation of ES systems. 
The simulation comparison between the ``water tank" and ``load/carbon-free generator" ES carbon footprint models are provided in Section \ref{sec:escomp}.

We note that
alongside the two proposed carbon footprint models for ES systems, there could be other valid models. The use of different carbon models can lead to discrepancies in carbon accounting results, operational decision-making, and market design for ES systems. Hence, it is crucial to assess the impact and select an appropriate model that aligns with practical objectives and specific requirements.

%% file: Simulation.tex
\section{Numerical Experiments} \label{sec:simulation}

In this section, we build a carbon-aware economic dispatch model based on the C-OPF method as an example for numerical tests, comparing it with conventional OPF-based schemes.

\subsection{C-OPF-Based and OPF-Based Economic Dispatch Models}\label{sec:ed}

Based on the C-OPF method \eqref{eq:copf}, we develop a carbon-aware economic dispatch (C-ED) model \eqref{eq:ed} as a specific application example. 
It
involves the multi-period optimal power scheduling of various generators and ES systems, while considering a complete AC power flow model, network operational constraints, and carbon flow equations and constraints. The C-ED model \eqref{eq:ed} 
is formulated as a nonconvex optimization problem and we solve it using the solver IPOPT \cite{biegler2009large}.
\begin{subequations}\label{eq:ed}
    \begin{align}
        \text{Obj.} &\  \min \sum_{i\in\sN} \sum_{t\in\sT}\Big[\sum_{g\in\sG_i} \!  \big(c_{i,g}^2  (P_{i,g,t}^{\rG})^2 + c_{i,g}^1  {P_{i,g,t}^{\rG}} + c_{i,g}^0\big)\hspace{-90pt} \nonumber\\
        & \qquad\qquad\qquad   + c_i^{\mathrm{es}} (P_{i,t}^{\mathrm{dc}} + P_{i,t}^{\mathrm{ch}}) 
 \Big],  \label{eq:ed:obj}\\
        \text{s.t.}\, 
        & \underline{P}_{i,g,t}^{\rG} \!\leq\! {P}_{i,g,t}^{\rG} \!\leq\!
        \bar{P}_{i,g,t}^{\rG}, && \hspace{-98pt} \forall i\!\in\! \sN, g\!\in\!\sG_i,t\!\in\!\sT\label{eq:ed:pgencap}
        \\
        & \underline{Q}_{i,g,t}^{\rG} \leq {Q}_{i,g,t}^{\rG} \leq \bar{Q}_{i,g,t}^{\rG}, &&\hspace{-98pt} \forall i\!\in\! \sN, g\!\in\!\sG_i,t\!\in\!\sT\label{eq:ed:qgencap}
        \\
        & \underline{\Delta}_{i,g}^{\rG} \!\leq\! {P}_{i,g,t}^{\rG} \!-\! {P}_{i,g,t\!-\!1}^{\rG} \!\leq\! \bar{\Delta}_{i,g}^{\rG}, && \hspace{-98pt} \forall i\!\in\! \sN, g\!\in\!\sG_i, t\!\in\!\sT\label{eq:ed:ramp} \\
         &  \hat{P}_{ij,t}^i - \hat{P}_{ji,t}^i     =  \big(V_{i,t}^2 -V_{i,t}V_{j,t}\cos(\theta_{i,t}-\theta_{j,t})\big)g_{ij} \nonumber\\
         &\quad  -V_{i,t}V_{j,t}\sin(\theta_{i,t}-\theta_{j,t})b_{ij},&& \hspace{-72pt} \forall ij\!\in\! \sE,t\!\in\!\sT   \label{eq:ed:Pfr}\\
          &  \hat{P}_{ij,t}^j - \hat{P}_{ji,t}^j     = - \big(V_{j,t}^2 -V_{j,t}V_{i,t}\cos(\theta_{j,t}-\theta_{i,t})\big)g_{ij} \nonumber\\
         &\quad+V_{j,t}V_{i,t}\sin(\theta_{j,t}-\theta_{i,t})b_{ij},&& \hspace{-72pt} \forall ij\!\in\! \sE,t\!\in\!\sT   \label{eq:ed:Pto}\\
 & Q_{ij,t}^i    =   \big(V_{i,t}V_{j,t}\cos(\theta_{i,t}-\theta_{j,t})-V_{i,t}^2\big)b_{ij}\nonumber \\
     &\quad  -V_{i,t}V_{j,t} \sin(\theta_{i,t}-\theta_{j,t})g_{ij},&& \hspace{-72pt} \forall ij\!\in\! \sE,t\!\in\!\sT \label{eq:ed:Qfr}\\
  & Q_{ij,t}^j    =   -\big(V_{j,t}V_{i,t}\cos(\theta_{j,t}-\theta_{i,t})-V_{j,t}^2\big)b_{ij}\nonumber \\
         &\quad  + V_{j,t}V_{i,t} \sin(\theta_{j,t}-\theta_{i,t})g_{ij},&& \hspace{-72pt} \forall ij\!\in\! \sE,t\!\in\!\sT \label{eq:ed:Qto} \\
        & \sum_{g\in \sG_i}\! P_{i,g,t}^{\rG}  \!-\!   P_{i,t}^{\rL} +P_{i,t}^{\mathrm{dc}} - P_{i,t}^{\mathrm{ch}} \nonumber\\
        &\qquad\qquad =\sum_{j\in \sN_i} (\hat{P}^i_{ij,t} \!-\!  \hat{P}^i_{ji,t}), &&\hspace{-72pt} \forall i\!\in\! \sN, t\!\in\!\sT \label{eq:ed:Pflow}\\
         & \sum_{g\in \sG_i} Q_{i,g,t}^{\rG}   -   Q_{i,t}^{\rL} =\sum_{j\in \sN_i} {Q}^i_{ij,t},   && \hspace{-72pt}  \forall i\!\in\! \sN, t\!\in\!\sT \label{eq:ed:Qflow}\\
        & (\hat{P}_{ij,t}^i - \hat{P}_{ji,t}^i )^2 + (Q_{ij,t}^i)^2\leq \Bar{S}_{ij}^2,&& \hspace{-72pt}  \forall ij\!\in\! \sE, t\!\in\!\sT\label{eq:ed:frthermal} \\
                & (\hat{P}_{ij,t}^j - \hat{P}_{ji,t}^j )^2 + (Q_{ij,t}^j)^2\leq \Bar{S}_{ij}^2,&& \hspace{-72pt}  \forall ij\!\in\! \sE, t\!\in\!\sT\label{eq:ed:tothermal} \\
           &\underline{V}_i \leq  V_{i,t} \leq \Bar{V}_i,\ \underline{\theta}_i \leq  \theta_{i,t} \leq \Bar{\theta}_i,&&\hspace{-72pt}  \forall i\!\in\!\sN, t\!\in\!\sT  \label{eq:ed:v}\\
        &0\!\leq\! \hat{P}_{ij,t}^i, 0\!\leq\!\hat{P}_{ji,t}^i, \, \hat{P}_{ij,t}^i\cdot\hat{P}_{ji,t}^i\!\leq\! 0, &&\hspace{-72pt}  \forall ij\!\in\! \sE, t\!\in\!\sT \label{eq:ed:frcomple}\\
      &0\!\leq\!\hat{P}_{ij,t}^j, 0\!\leq\! \hat{P}_{ji,t}^j,\, \hat{P}_{ij,t}^j\cdot\hat{P}_{ji,t}^j\!\leq\! 0, &&\hspace{-72pt}  \forall ij\!\in\! \sE, t\!\in\!\sT  \label{eq:ed:tocomple}\\
       &   w_{i,t} \leq \bar{w}_{i,t},   && \hspace{-72pt} \forall i\!\in\! \sN,   t\!\in\!\sT\label{eq:ed:cfc} \\
       &\text{Equations \eqref{eq:ESpower}, \eqref{eq:ES:w}, \eqref{eq:cf:ES1}}.\label{eq:ed:ES}
    \end{align}
\end{subequations}
The objective \eqref{eq:ed:obj} aims to minimize the total generation and ES operational costs, where $c_i^{\mathrm{es}}$ denotes the ES degradation cost coefficient.
Equations \eqref{eq:ed:pgencap}-\eqref{eq:ed:ramp} denote the active and reactive power generation capacity limits and the ramping constraints for various generators.  Equations \eqref{eq:ed:Pfr} and \eqref{eq:ed:Pto} represent the full AC power flow equations for the active power flow values at the sending and receiving nodes, respectively. Here, we employ the dual power flow reformulation method introduced in Section \ref{sec:keyissue} to address the issue of unknown power flow directions. As shown in Figure \ref{fig:cfvalue}, ${P}_{ij,t}^i\!=\!\hat{P}_{ij,t}^i \!-\! \hat{P}_{ji,t}^i$ and ${P}_{ij,t}^j\!=\!\hat{P}_{ij,t}^j \!-\! \hat{P}_{ji,t}^j$ are the active power flow values of branch $ij$ that are measured at node $i$ and node $j$, respectively. Thus, the power loss of branch $ij$ is $P_{ij,t}^{\mathrm{loss}}\!=\! |{P}_{ij,t}^i-{P}_{ij,t}^j|$. Similarly,
equations \eqref{eq:ed:Qfr} and \eqref{eq:ed:Qto} represent the full AC power flow equations for the branch reactive power flows at the sending and receiving nodes, respectively. Equations \eqref{eq:ed:Pflow} and \eqref{eq:ed:Qflow} are the active and reactive power balance constraints at each node. Equations \eqref{eq:ed:frthermal} and \eqref{eq:ed:tothermal} represent the line thermal constraints at the sending and receiving nodes, respectively. Equation \eqref{eq:ed:v} denotes the upper and lower limits on nodal voltage magnitudes and phase angles. Equations \eqref{eq:ed:frcomple} and \eqref{eq:ed:tocomple} represent the nonnegativity and complementarity constraints for the dual power flow values; see Section \ref{sec:keyissue} for more explanations. 
 Equation
\eqref{eq:ed:cfc} is the carbon flow constraint that imposes a cap on the nodal carbon intensities. This constraint enables the active management of nodal carbon intensities, ensuring that clean power is supplied to users with a carbon intensity no greater than $\bar{w}_{i,t}$. 
Equation \eqref{eq:ed:ES} collects
the dynamic power model \eqref{eq:ESpower}  and the ``water tank" carbon emission model \eqref{eq:ES:w} for ES systems, as well as 
the carbon flow equations \eqref{eq:cf:ES1}.

For comparison, we also build a conventional OPF-based ED model that does not incorporate the carbon flow model and constraints. Essentially, the OPF-based ED model is a reduced version of the C-ED model \eqref{eq:ed}, which excludes the carbon-related constraints \eqref{eq:ed:cfc}, \eqref{eq:ES:w}, \eqref{eq:cf:ES1} and does not need the introduction of dual power flow variables. The OPF-based ED model is also a nonconvex optimization problem due to the full AC power flow equations and the dynamic ES power model, and we solve it using the solver IPOPT \cite{biegler2009large}.

\subsection{Test System and Simulation Settings}

\begin{figure}
    \centering
    \includegraphics[scale=0.38]{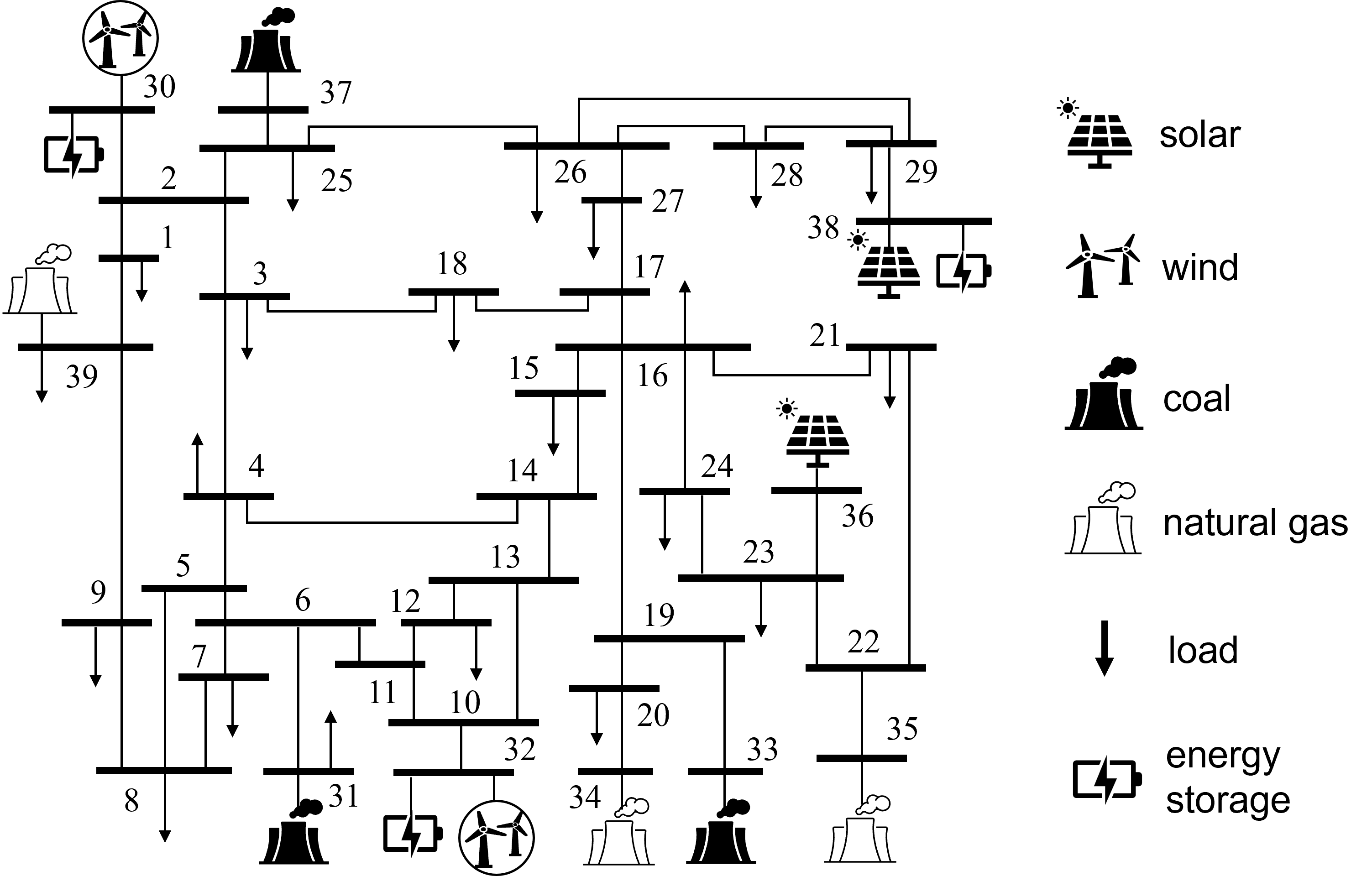}
    \caption{The modified New England 39-bus test system.}
    \label{fig:39bus}
\end{figure}

In the simulations, we consider day-ahead economic dispatch with $T\!=\!12$ time steps and 2-hour time intervals. The modified New England 39-bus  system, as shown in Figure \ref{fig:39bus}, is used as the test system, which includes 3 coal power plants, 3 natural gas power plants, 2 wind farms, 2 solar farms, 3 ES systems, and 21 loads.
The generation carbon emission factors $w_{i,g}^{\mathrm{G}}$ are 2.26, 0.97, and 0 (lbs/kWh) for coal plants, natural gas plants, and renewable generators, respectively. 
The carbon flow constraint \eqref{eq:ed:cfc} is imposed for all the load nodes and all $t\in \sT$, and we set the cap $\bar{w}_{i,t} = 1.2$ lbs/kWh.\footnote{In the simulations, we set the same carbon intensity cap for all load nodes and all times for simplicity.
In practical applications, distinct carbon intensity caps can be implemented for different nodes to distinguish the ``cleanness" level of electricity supply at specific locations.} For each ES system, we set the storage efficiency factor as $\kappa_i=0.99$ and the charging and discharging efficiency coefficients as $\eta_i^{\mathrm{ch}} = \eta_i^{\mathrm{dc}}=0.98$. 
Other 
detailed system parameters and settings are provided in \cite{ieee39}. 

\subsection{Simulation Results Comparison of C-OPF and OPF}

\subsubsection{Nodal Carbon Intensity} 
\begin{figure}
    \centering
            \includegraphics[scale=0.335]{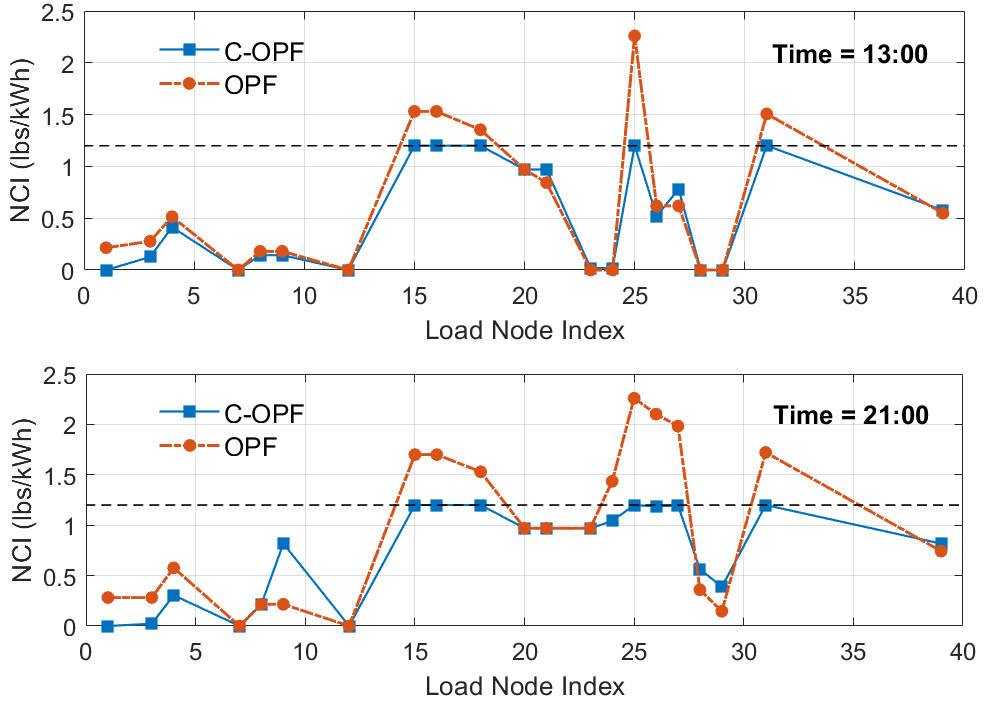}
    \caption{The nodal carbon intensities (NCI) for all load nodes under the C-OPF-based and OPF-based economic dispatch schemes at 13:00 and 21:00. (The black dashed line denotes the NCI cap $\bar{w}_{i,t}\!=\!1.2$ lbs/kWh).}
    \label{fig:nci}
\end{figure}

\begin{figure}
    \centering
    \includegraphics[scale=0.5]{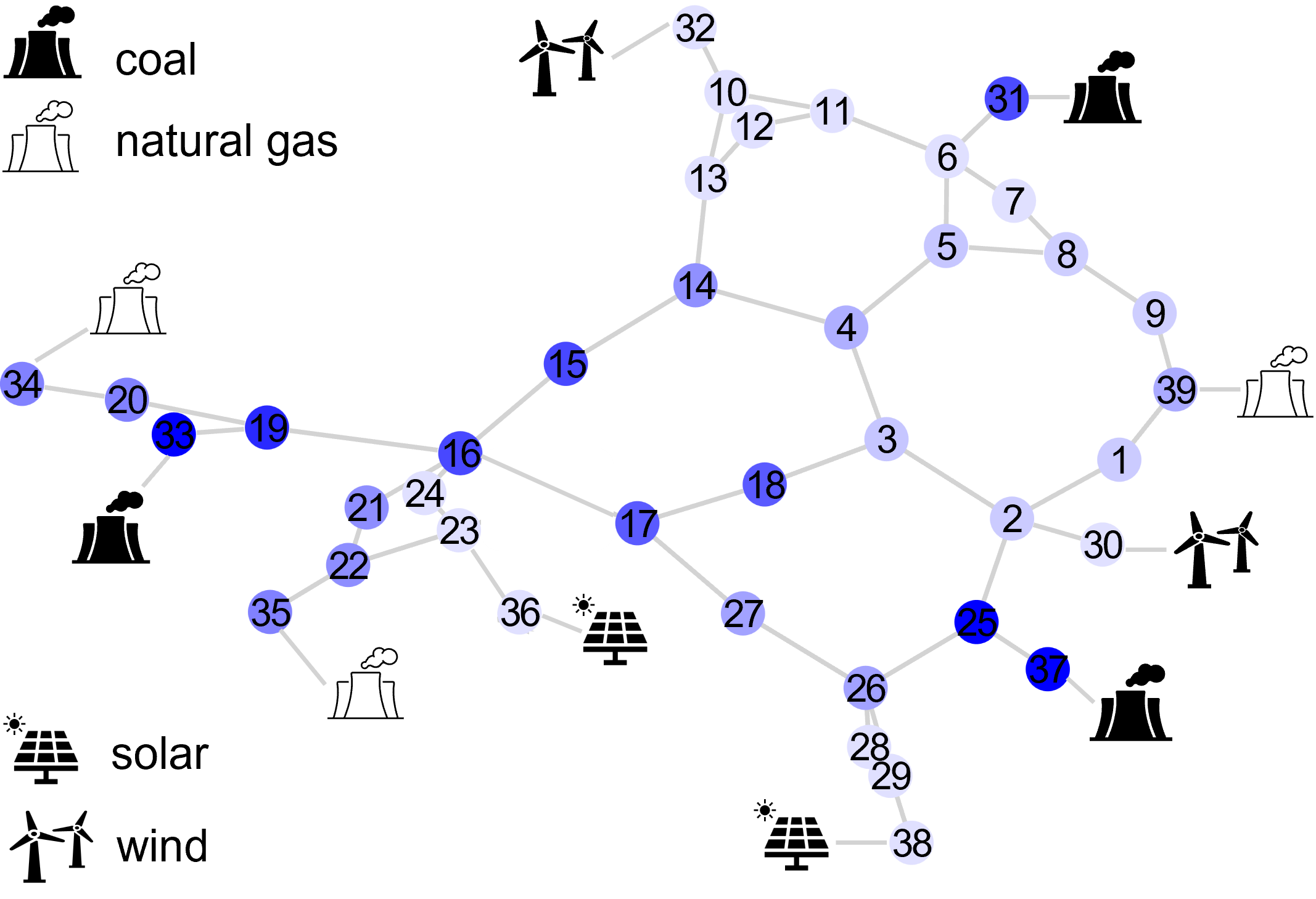}
    \caption{Visualization of nodal carbon intensities of the modified New England 39-bus test system at 13:00. (Darker color indicates higher carbon intensity).}
    \label{fig:gridnci}
\end{figure}

Figure \ref{fig:nci} illustrates the nodal carbon intensities for all load nodes under the C-OPF-based and OPF-based ED schemes, respectively. It shows that the C-OPF-based ED model can generate effective power dispatch schemes that keep the nodal carbon intensities of load nodes below the cap of $\bar{w}_{i,t} \!=\! 1.2$ lbs/kWh. In contrast, the OPF-based ED scheme frequently exceeds the nodal carbon intensity cap. Moreover, the grid's nodal carbon intensities exhibit significant temporal variation and spatial diversity.  From  Figure \ref{fig:nci}, 
it is observed that the nodal carbon intensities at 13:00 are generally lower than those at 21:00, due to the higher penetration of renewable generation at 13:00, as illustrated in Figure \ref{fig:generation}. Figure \ref{fig:gridnci} visualizes the grid's nodal carbon intensities at 13:00, with darker blue colors indicating higher carbon intensity at each node. It demonstrates that the nodal carbon intensities calculated using the carbon flow method can reflect the proximity to different fuel types of generation and align with the physical power flow. In comparison, the grid average emission factor at 13:00 is calculated to be 0.5 lbs/kWh, which only gives an indication of the overall grid average emission state and falls short of providing detailed insight into local emissions at different locations.

\subsubsection{Power Dispatch Scheduling Decisions}
\begin{figure}
     \centering
     \includegraphics[scale=0.335]{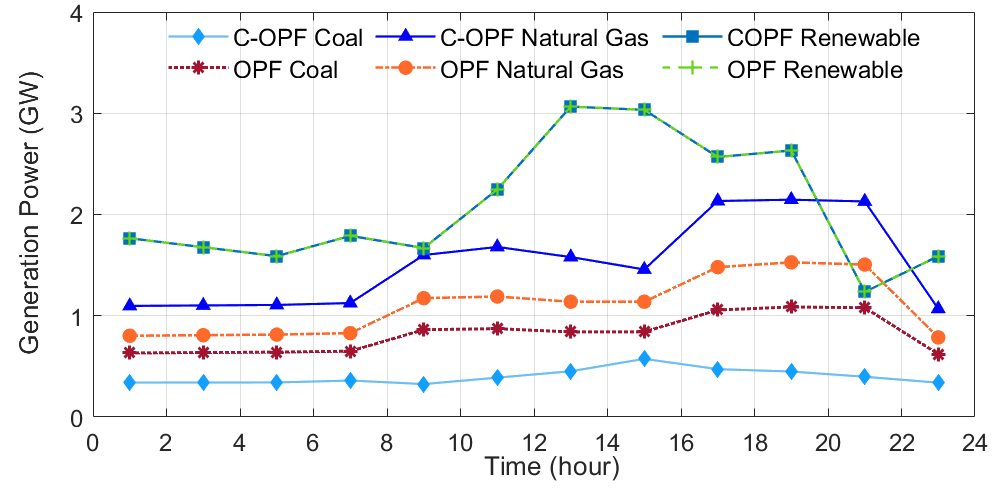}
     \caption{The generation outputs of various power plants over 24 hours under the C-OPF-based and OPF-based economic dispatch schemes.}
     \label{fig:generation}
 \end{figure}

Figure  \ref{fig:generation} illustrates the generation decisions of the C-OPF-based and OPF-based ED schemes over a 24-hour period. In both schemes, the renewable generation, i.e., solar and wind generation, is fully utilized without curtailment. The primary distinction is that the C-OPF-based ED scheme results in more generation from expensive yet clean natural gas plants and less generation from cheap but high-emission coal plants to meet the carbon emission constraints, compared with the OPF scheme. The operational costs of these two ED schemes are presented in Table \ref{tab:cost}. It is seen that the total operational cost of the C-OPF-based ED scheme is higher than that of the OPF scheme, due to the increased generation from natural gas plants as a substitute for coal plant generation. 

 \begin{table}[]
 \caption{Operational cost of OPF-based and C-OPF-based ED schemes. }
\begin{tabular}{cccccc}
\hline
 Cost (k\$)     & Coal   & Natural Gas & Renewable & \begin{tabular}[c]{@{}c@{}}Energy \\ Storage\end{tabular} & Total \\ \hline \\[-0.8em]
OPF                                                            &4515.4 & 6734.1      & 50.7      & 15.1      & 11315.3      \\ \\[-0.8em]
C-OPF                                                         & 1345.0 & 12116.3     & 50.7      & 15.7      & 13527.7    \\ \hline
\end{tabular}
\label{tab:cost}
\end{table}

\subsubsection{Grid Carbon Emissions and Energy Storage}
\begin{figure}
    \centering
    \includegraphics[scale=0.335]{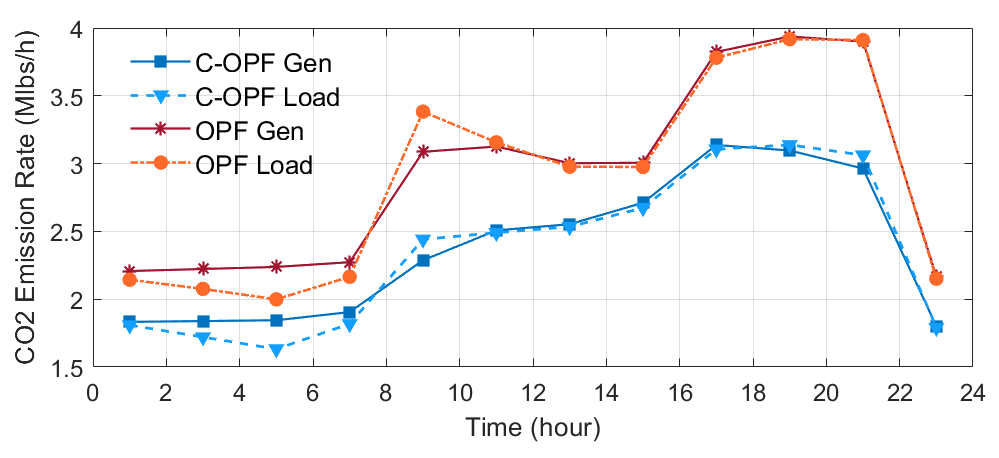}
        \vspace{-10pt}
    \caption{The total (Scope 1) generation-side carbon emission rate and the total (Scope 2) attributed load-side carbon emission rate over time under the C-OPF-based and OPF-based economic dispatch schemes.}
    \label{fig:carbonemi2}
\end{figure}

Figure \ref{fig:carbonemi2} illustrates the total generation-side and attributed demand-side carbon emission rates over time.  It is observed that 
the C-OPF-based ED scheme yields reduced carbon emission rates on both the generation side and demand side, compared with the OPF-based ED scheme. 
That is because the carbon flow constraint \eqref{eq:ed:cfc} in the C-OPF model requires a larger share of power supply from clean generators to meet the nodal carbon intensity cap. As a result, the total system emissions over 24 hours are reduced from the amount of 69,981.9 klbs in the OPF scheme to 56,958.8 klbs in the C-OPF scheme.

\begin{figure}
     \centering
\includegraphics[scale=0.334]{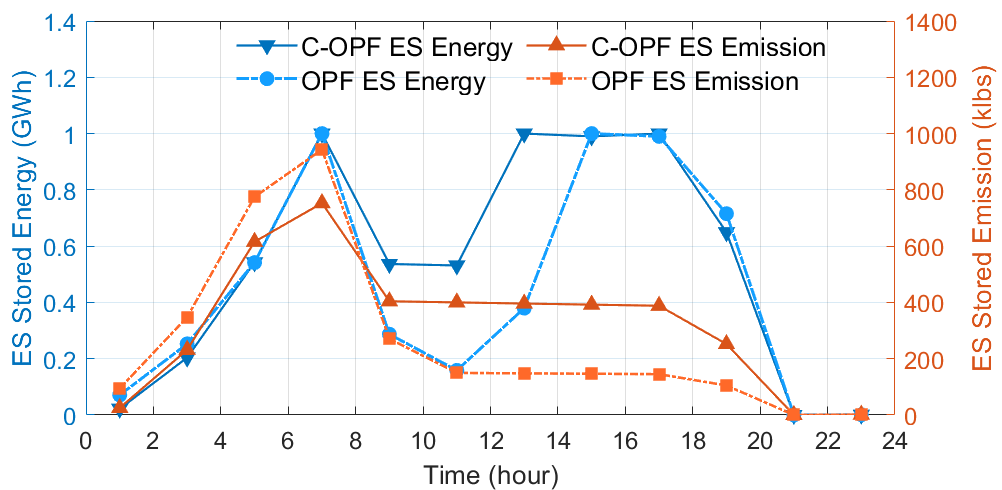}
     \vspace{-10pt}
     \caption{The stored energy and virtual carbon emissions of the ES system at node-38 under the C-OPF-based and OPF-based economic dispatch schemes using the ``water tank" ES carbon footprint model. }
     \label{fig:ESeE}
 \end{figure}

Moreover, slight differences between the generation-side emission rates and the attributed load-side emission rates are observed in Figure \ref{fig:carbonemi2} for both the OPF and C-OPF schemes. These differences result from the attributed emission rates for the operation of ES systems and power loss. 
Consistent with the power conservation law, 
our proposed carbon accounting mechanisms based on the carbon flow method ensure the ``\emph{carbon conservation principle}" \cite{ourvisionpaper} that the total generation-side carbon emissions equal the sum of emissions attributed to the total load, power loss, and ES systems at all times.  
Figure \ref{fig:ESeE} illustrates the time trajectories of stored energy $e_{i,t}^{\mathrm{es}}$ and virtual carbon emissions $E_{i,t}^{\mathrm{es}}$ of the ES system at node-38. 
It is observed that the virtual carbon emissions generally increase when the ES system charges and decrease when it discharges. The period from 11:00 to 15:00, when charging does not result in increased emissions, occurs because the ES system charges with carbon-free renewable electricity.
By combining Figures \ref{fig:carbonemi2} and \ref{fig:ESeE}, it is seen that from 1:00 to 7:00, the generation-side emissions are higher than the load-side emissions due to the ES systems charging and absorbing emissions from the grid. Conversely, at 9:00, the ES systems are discharging and injecting emissions back into the grid, resulting in higher load-side emissions than the generation-side emissions.

Table \ref{tab:emi} summarizes the carbon accounting results across 24 hours for the C-OPF-based and OPF-based ED schemes, where the (Scope 2) attributed emissions are calculated using the mechanisms introduced in Sections  \ref{sec:accounting} and \ref{sec:watertank}. 
It is verified that the total 
 (Scope 1) generation-side emissions are equal to the total (Scope 2) emissions attributed to loads, grid power loss, and ES systems.

 \begin{table}[H]
 \caption{Carbon accounting results based on carbon flow method and ``water tank" ES carbon footprint model.}
\begin{tabular}{lcccc}
\hline \\[-0.8em]
\multirow{2}{*}{(klbs)} & \multirow{2}{*}{\begin{tabular}[c]{@{}c@{}}(Scope 1) Generation\\ -Side Emissions\end{tabular}}  & \multicolumn{3}{c}{(Scope 2) Attributed Emissions} \\  \cline{3-5} \\[-0.8em] &                          & Load        & Power Loss        & ES Systems       \\ \\[-0.8em] \hline \\[-0.8em]
OPF & 69,982    & 69,273           & 639                 & 70                \\ \\[-0.8em]
C-OPF    & 56,959  & 56,433   & 457                 & 69                \\ \hline
\end{tabular}
\label{tab:emi}
\end{table}

\subsection{Impact of Nodal Carbon Intensity Cap}

\begin{figure}
    \centering
    \includegraphics[scale=0.335]{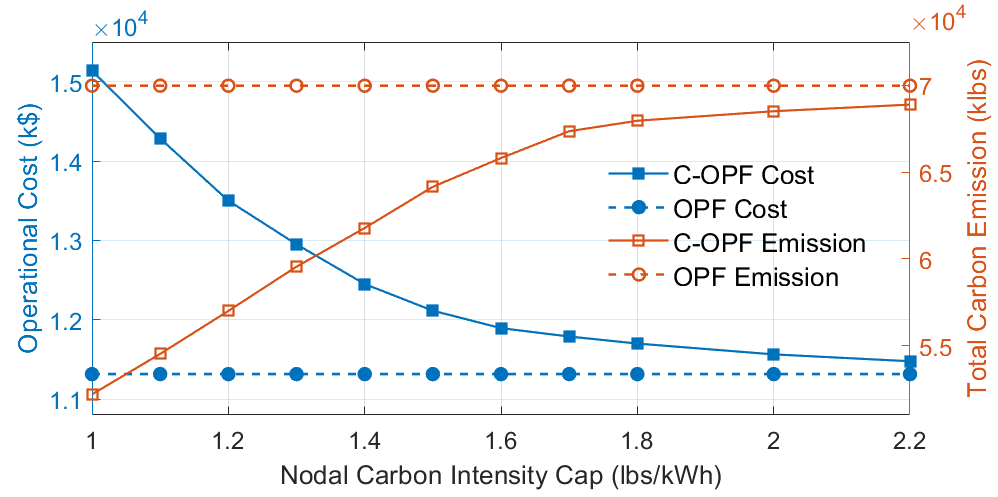}
    \caption{Total system emissions and total operational costs of OPF-based and C-OPF-based ED schemes under different nodal carbon intensity caps $\bar{w}_{i,t}$.}
    \label{fig:cap}
\end{figure}

To study the impact of the carbon flow constraint \eqref{eq:ed:cfc}, we adjust the nodal carbon intensity cap $\bar{w}_{i,t}$ from $1$ lbs/kWh to $2.2$ lbs/kWh uniformly for all load nodes, and run the C-OPF-based ED model for each case. Figure \ref{fig:cap} illustrates the results of total system emissions and total operational costs under different carbon intensity caps. 
 It is observed that as the nodal carbon intensity cap increases, the total system operational cost of the C-OPF scheme decreases monotonically and converges to that of the OPF scheme; simultaneously, the total system carbon emissions increase and also converge to those of the OPF scheme. Essentially, the C-OPF-based ED scheme strikes an optimal trade-off between operational cost and carbon emission reduction. 
This is also consistent with the intuition that C-OPF solutions gradually reduce to traditional OPF solutions as carbon flow constraints become less restrictive.

\subsection{Comparison of Two ES Carbon Footprint Models}\label{sec:escomp}

In the simulations above, the ``water tank" ES carbon footprint model is adopted in the C-OPF-based ED model and the carbon flow calculation for demand-side carbon accounting. In contrast, this subsection implements the ``load/carbon-free generator (LCG)" ES carbon footprint model for comparison. The carbon accounting results using the LCG ES carbon footprint model are shown in Table \ref{tab:emi2}. 

\begin{table}[H]
 \caption{Carbon accounting results based on carbon flow method and ``load/carbon-free generator" ES carbon footprint model.}
\begin{tabular}{lcccc}
\hline \\[-0.8em]
\multirow{2}{*}{(klbs)} & \multirow{2}{*}{\begin{tabular}[c]{@{}c@{}}(Scope 1) Generation\\ -Side Emissions\end{tabular}}  & \multicolumn{3}{c}{(Scope 2) Attributed Emissions} \\  \cline{3-5} \\[-0.8em] &                          & Load        & Power Loss        & ES Systems       \\ \\[-0.8em] \hline \\[-0.8em]
OPF &    69,982  &     68,375        &            631       &    976              \\ \\[-0.8em]
C-OPF    & 57,322  &  56,081   &                  463   &    778          \\ \hline
\end{tabular}
\label{tab:emi2}
\end{table}

Compared with Table \ref{tab:emi} that uses the ``water tank" model, it is seen from Table \ref{tab:emi2} that for the OPF model, its optimal solution and generation-side emissions remain the same, as the OPF model does not involve carbon emissions. However, the attributed emissions calculated using the carbon flow method are different. Specifically, in the LCG ES carbon model, a greater share of the carbon footprints is attributed to the ES systems, thereby reducing the carbon footprints attributed to the loads and power loss. This is because, unlike the ``water tank" model, ES systems are treated as pure loads during charging in the LCG model. As a result, virtual carbon emissions absorbed from the grid accumulate locally at the ES systems and are not released back to the grid during discharging, which become the carbon footprints of the ES systems. In terms of the C-OPF model, using the LCG ES carbon footprint model alters its optimal solution, leading to higher generation-side emissions and a reduced operational cost of $13,506.2$ k\$. Because the ES systems impact the grid's carbon flow only during discharging, when the LCG ES carbon footprint model treats ES systems as carbon-free generators that lower the grid's carbon intensities. Thus, it allows for increased generation from high-emission but cheaper coal plants, while still satisfying the nodal carbon intensity constraints. Similarly, in the C-OPF case, much higher carbon footprints of 778 klbs are attributed to the ES systems, in contrast to the 69 klbs shown in Table \ref{tab:emi} under the ``water tank" model.

\subsection{Computational Efficiency}

The numerical experiments are implemented in a computing environment with Intel(R) Core(TM) i7-1185G7 CPUs running at 3.00 GHz and with 16 GB RAM. We use the JuMP language \cite{Lubin2023} in Julia to build optimization models and solve them using the IPOPT solver (version 3.14.10) \cite{biegler2009large}. It takes 176.6 seconds on average to solve the C-OPF-based ED model and 55.6 seconds to solve the OPF-based ED model. 

 Additionally, we implement the DC power flow model \eqref{eq:dc} instead of the full AC power flow model in the C-OPF-based and OPF-based ED models for simulation comparison. In this case, both ED models are significantly simplified due to omitting voltage magnitude and reactive power, neglecting power line losses, and replacing nonlinear power flow equations with linear DC flow equations. Under the DC power flow model, 
 it takes 16.4 seconds on average to solve the C-OPF-based ED model and 1.2 seconds to solve the OPF-based ED model.

Due to the dual power flow reformulation and the addition of carbon flow equations and constraints, the C-OPF model generally has a larger problem size and requires more solution time than the traditional OPF model. In practice, several methods can improve the solution efficiency of the C-OPF model. For instance, by inspecting all branches and identifying those whose power flow directions can be predetermined based on the network configuration, the dual power flow reformulation can be avoided for these branches. Additionally, the solutions of the OPF model can be employed to warm-start the C-OPF model. 
A key future research direction is to develop efficient linearization and convexification approaches for the nonconvex carbon flow equations to fundamentally enhance the solution efficiency of the C-OPF model.

%% file: Conclusion.tex
\section{Conclusion}\label{sec:conclusion}

This paper proposes a generic Carbon-aware Optimal Power Flow (C-OPF) methodology as a
fundamental tool for guiding decarbonization decision-making in electric power systems.  
As a carbon-aware generalization of conventional OPF models, the C-OPF model 
enables the joint management of the grid's carbon footprints and power flows. A reformulation technique is introduced to address the issue of unknown power flow directions in the C-OPF model. Additionally, we propose two novel carbon footprint models for ES systems as well as their corresponding 
carbon accounting mechanisms, facilitating optimal carbon-aware ES operation. Numerical
simulations demonstrate that C-OPF-based schemes can effectively coordinate diverse energy resources for grid decarbonization while ensuring that decisions comply with regulatory requirements on carbon footprints. 
This paper represents preliminary work introducing the emerging and promising technique of C-OPF for supporting optimal power system decarbonization decisions.   Extensive future research is anticipated to further advance the methodology of C-OPF. Potential future directions include 1) theoretical advances in the C-OPF modeling and optimization, such as efficient linearization and convexification of the C-OPF model, and 2) practical applications of C-OPF to power grid decision-making problems, such as carbon-aware demand response and carbon-electricity pricing.

%% file: Appendix.tex
\appendices

\section{Derivation of Carbon Flow Matrix Form \eqref{eq:mci:ma}} \label{app:cfder}

From \eqref{eq:mci}, we first multiply both sides by the denominator of \eqref{eq:mci} and move the term $\sum_{k\in \sN_i^+ } w_k P_{ki}^i$ to the left-hand side, leading to \eqref{eq:der:step1} for all $i\in\sN$:
\begin{align} \label{eq:der:step1}
  w_i \Big(\underbrace{\sum_{g\in\sG_i}  P_{i,g}^{\rG} + \sum_{k\in \sN_i^+ }\!\! P_{ki}^i}_{:= P_i^{\mathrm{in}}}\Big) -  \sum_{k\in \sN_i^+ }\!\! w_k P_{ki}^i =   R_i^{\rG}, 
\end{align}
where $R_i^{\rG}\!:=\! \sum_{g\in\sG_i}\! w_{i,g}^{\rG} P_{i,g}^{\rG}$ denotes the total generation carbon emission rate at node $i$. We then stack up equation \eqref{eq:der:step1} for all nodes $i\in\sN$ into a column form. On the left-hand side, the first term of \eqref{eq:der:step1} becomes $\bm{P}_{\rN}\bm{w}_{\rN}$, where $\bm{w}_{\rN}\!:=\! (w_i)_{i\in\sN}$ is the column vector that collects the nodal carbon intensities $w_i$, and $\bm{P}_{\rN}\!:=\!\text{diag}(P_i^{\mathrm{in}})$ is the diagonal matrix whose $i$-th diagonal entry is the nodal active power inflow $P_i^{\mathrm{in}}$ to node $i$. The second term of \eqref{eq:der:step1} becomes $-\bm{P}_{\rB}\bm{w}_{\rN}$, where $\bm{P}_{\rB} \!\in\! \R^{N\times N}$ is the branch power inflow matrix that is constructed by letting $\bm{P}_{\rB}[i,k] \!=\! P_{ki}^i$ and $\bm{P}_{\rB}[k,i] \!= \!0$ if node $k$ sends power flow $P_{ki}^i$ to node $i$. The right-hand side of \eqref{eq:der:step1} becomes the column vector $\bm{r}_{\rG}\!:=\! (R_i^{\rG})_{i\in\sN}$. Thus, equation \eqref{eq:der:step1} for all $i\in\sN$ can be equivalently reformulated as \eqref{eq:der:step2}:
\begin{align}\label{eq:der:step2}
  (\bm{P}_{\rN} -\bm{P}_{\rB} )\bm{w}_{\rN} = \bm{r}_{\rG}.
\end{align}
Taking the matrix inverse of \eqref{eq:der:step2} leads to \eqref{eq:mci:ma}. See \cite{kang2015carbon} for more derivation details.

%% file: Reference.tex
\bibliography{IEEEabrv, ref}

%% file: BioNoPhoto.tex
\vskip -4pt plus -1fil

\begin{IEEEbiographynophoto}{Xin Chen} is an Assistant Professor in the Department of Electrical and Computer Engineering at Texas A\&M University (TAMU). Dr. Chen directs the Smart Power, Energy and Decision-making (SPEED) Lab at TAMU. His research lies in the intersection of control, machine/reinforcement learning, and optimization for smart sustainable power and energy systems. He received the Ph.D. degree in electrical engineering from Harvard University, the master’s degree in electrical engineering and two bachelor’s degrees in engineering and economics from Tsinghua University. He was a Postdoctoral Associate affiliated with MIT Energy Initiative at Massachusetts Institute of Technology. He is a recipient of the IEEE PES Outstanding Doctoral Dissertation Award, IEEE Transactions on Smart Grid Top-5 Outstanding Papers, and multiple best paper awards at IEEE control and power conferences.  
\end{IEEEbiographynophoto}

\vskip 0pt plus -1fil

\begin{IEEEbiographynophoto}{Andy Sun} is is the Iberdrola-Avangrid Professor in Electric Power Systems in the Sloan School of Management at the Massachusetts Institute of Technology (MIT). He received his PhD degree in Operations Research from MIT and was a postdoctoral fellow in the Mathematical Sciences division of the IBM T. J. Watson Research Center at Yorktown Heights, NY. He was a faculty member in the School of Industrial and Systems Engineering at Georgia Institute of Technology before joining MIT. Dr. Sun's research focuses on optimization and computation for large-scale electric power system control, operations, and planning, market design for power grid decarbonization, and electrification of transportation. His research on robust operation of power grids has been influential in improving reliability unit commitment in ISO/RTO markets.  
\end{IEEEbiographynophoto}

\vskip 0pt plus -1fil

\begin{IEEEbiographynophoto}{Wenbo Shi} is the Founder/CEO of Singularity Energy, a startup that offers advanced carbon and clean energy management software and data solutions for utilities, grid operators,
corporations and technology providers to accurately measure emissions and optimize their decision-making for grid decarbonization. Singularity is proud to partner with organizations such as Southern Company, MISO, Eversource Energy, Enersponse, Sense, and Measurable.Before founding Singularity, Wenbo was a postdoctoral researcher at Harvard University. He received his Ph.D. from University of California, Los Angeles in 2015.
\end{IEEEbiographynophoto}

\vskip 0pt plus -1fil

\begin{IEEEbiographynophoto}{Na Li}
 is a Winokur Family Professor of Electrical Engineering and Applied Mathematics at Harvard University.  She received her Bachelor's degree in Mathematics from Zhejiang University in 2007 and Ph.D. degree in Control and Dynamical systems from California Institute of Technology in 2013. She was a postdoctoral associate at the Massachusetts Institute of Technology 2013-2014.  She has held a variety of short-term visiting appointments including the Simons Institute for the Theory of Computing, MIT, Google Brain, and MERL. Her research lies in the control, learning, and optimization of networked systems, including theory development, algorithm design, and applications to real-world cyber-physical societal system. 
\end{IEEEbiographynophoto}